\documentclass[12pt,a4paper]{article}
\usepackage{amsmath,amssymb,amsthm}

\newcommand{\re}{{\mathbb R}}
\newcommand{\n}{{\mathbb N}}

\newcommand{\z}{{\mathbb Z}}
\newcommand{\cA}{{\mathcal{A}}}

\newcommand{\cF}{{\mathcal{F}}}
\newcommand{\cG}{{\mathcal{G}}}

\newcommand{\cL}{{\mathcal{L}}}

\newcommand{\cT}{{\mathcal{T}}}

\newcommand{\cS}{{\mathcal{S}}}

\newcommand{\bA}{{\boldsymbol{A}}}

\newcommand{\bT}{{\boldsymbol{T}}}

\newcommand{\bc}{{\boldsymbol{c}}}

\newtheorem{theorem}{Theorem}
\newtheorem{prop}{Proposition}
\newtheorem{lemma}{Lemma}
\newtheorem{cor}{Corollary}
\newtheorem{remark}{Remark}
\newtheorem{ex}{Example}
\newtheorem{defi}{Definition}

\newtheorem{assum}{Assumption}

\date{}

\author{Maria Charina
\thanks{Fakult\"at f\"ur Mathematik,  Universit\"at Wien, Austria
{e-mail: \tt\small  maria.charina@univie.ac.at}}
and Vladimir Yu.~Protasov
\thanks{DISIM of University of L'Aquila and Department of Mechanics and Mathematics of Moscow State University  {e-mail: \tt\small
v-protassov@yandex.ru}}}

\title{Regularity of anisotropic refinable functions
\thanks{
The first author is sponsored by the Austrian Science Foundation (FWF) grant P28287-N35; the second author is supported by RFBR grants nos. 16-04-0832 and 17-01-00809.}}

\begin{document}
\maketitle

\begin{abstract}
This paper presents a detailed regularity analysis of multivariate refinable functions with general dilation matrices,
with emphasis on the anisotropic case. 
In the univariate setting, the smoothness of refinable functions is well understood by
means of the matrix approach.  In the multivariate setting, this approach has been extended only to
the special case of isotropic refinement with the dilation matrix all of whose eigenvalues are equal in the absolute value.  The general  anisotropic case has resisted to be fully understood:  the matrix approach can determine whether a refinable function belongs to $C(\re^s)$ or $L_p(\re^s)$, $1 \le p < \infty$,
but its H\"older regularity remained mysteriously unattainable.

We show how to compute the H\"older regularity in $C(\re^s)$ or $L_p(\re^s)$, $1 \le p < \infty$.  In the anisotropic case, our expression for the exact H\"older exponent of a refinable function reflects the impact of the variable moduli of the eigenvalues of the corresponding dilation matrix. In the isotropic case, our
results reduce to the well-known facts from the literature. We also analyze the higher regularity of anistropic refinable functions. We illustrate our results with several examples.

\smallskip

\smallskip

\begin{flushright}
\noindent{\bf Classification (MSCS): 65D17, 15A60, 39A99 }
\end{flushright}

\end{abstract}
\bigskip

\section{Introduction}\label{s.intr}

We study the multivariate refinement equation
\begin{equation}\label{eq.ref}
\varphi(x)\ = \ \sum_{k \in \z^s} \, c_k \, \varphi \, (Mx \, - \, k), \quad  x \in \re^s,
\end{equation}
with a compactly supported sequence $c=\{c_k \ : \ k \in \z\}$ of coefficients $c_k \in \re$ and with a general integer dilation
matrix~$M \in \z^{s \times s}$ all of whose eigenvalues are larger than one in the absolute value.
We do not make any assumptions on the stability of the integer shifts of $\varphi$.

In this paper, we characterize continuous and  $L_p$ solutions of \eqref{eq.ref}.
Our main contribution is the exact expression for the H\"older exponent of $\varphi$ in $C(\re^s)$ and in
$L_p(\re^s)$, $1 \le p < \infty$, see Theorems \ref{th.holder} and \ref{th.holder-p}. In the anisotropic case,
the H\"older exponent of  $\varphi$ reflects the influence of the invariant subspaces of $M$ corresponding to its different by modulus eigenvalues. In the isotropic case, when all the eigenvalues of $M$ are equal in the absolute value,
our results reduce to the well-known ones from the literature. We also analyze the higher regularity of continuous $\varphi$.

It is well known that compactly supported solutions ({\em refinable functions}) of \eqref{eq.ref}
can generate systems of multivariate wavelets or frames, see e.g~\cite{Chr, Chui, Daub, Prot_book, NPS}.  Refinable functions are
building blocks for the limits of subdivision algorithms widely used in approximation and for generating curves and
surfaces, see e.g.~\cite{CDM, ChuiV, DynLevin02, ReifPeter08, W}.
Refinable functions naturally appear in recent applications in combinatorics, multigrid, number theory, and probability
~\cite{CharinaDRT, CharinaDRT2, DDL, FS, KM, P00, P16}.

In the univariate case, $M\ge 2$ is an integer, there are several efficient methods for determining the
regularity of refinable functions. In \cite{CD, DD, E} the authors compute precisely the
Sobolev exponent of $\varphi \in L_2(\re)$.  The so-called {\em matrix approach} yields the H\"older exponent of $\varphi \in C(\re)$ and, in addition, provides a detailed
analysis of its moduli of continuity and of its local regularity~\cite{CH, DL,P06, Rio}. An obstacle to the
practical use of the matrix approach is the NP-hardness of the joint spectral radius computation.
This problem, however, was successfully resolved for a large class of problems by recent results in~\cite{GP1, GP2, MR}
where the authors presented fast and efficient methods of the joint spectral radius computation. Indeed, the invariant polytope algorithm~\cite{GP1} estimates the joint spectral radius for the corresponding transition matrices of size up to~$20$ and, in most cases, even
determines its precise value.

The generalization of the matrix approach to the multivariate case turned out to be a difficult task in the
case of general dilation matrices.  The special case of isotropic dilation is currently fully understood,
see~\cite{CDM, Charina, CJR, CD, DD, DynLevin02, E, HJ, H1, H2, J, JJ, JZ, JiangLiZhu, RonS}. Several partial results in the anisotropic case are also available:
for characterizations of continuity and $L_p$, $1 \le p < \infty$, regularity of  $\varphi$ see e.g.~\cite{CHM2, H3}; for estimates
for the H\"older exponent of $\varphi$ see e.g.~\cite{CHM2, H3}.

The reason for the difficulty of the anisotropic case is natural and hardly avoidable. In the univariate case, say $M = 2$,
the distance between two points $x, y\in \re$ can be expressed in terms of their binary expansions.
The distance between the values $\varphi(x)$ and $\varphi(y)$ depends on the behavior of the products of certain square matrices
derived from $c_k$, $k \in\z$. These two observations establish a correlation between $|x-y|$ and $|\varphi(x) - \varphi(y)|$,
which leads to the formula for the H\"older exponent of $\varphi$. This summarizes the essence of the matrix approach.
In the multivariate case, one can similarly estimate the distance between $\varphi(x)$ and $\varphi(y)$
by suitable matrix products. The problem occurs at an unexpected point:
the expression for the distance between $x, y \in \re^s$. One can try to use the corresponding $M$-adic expansions
with a certain set of digits from~$\z^s$, but such expansions do not provide a clear estimate for the distance between $x$ and $y$.
Indeed, unless the matrix $M$ is isotropic, multiplication by a high power of $M$ can enlarge  distances differently in different
directions. Hence, the points $M^\ell x$ and $M^\ell y$, $\ell \in \n$, whose $M$-adic expansions are essentially the same,
may have different asymptotic behavior as $\ell \to \infty$. Remarkably simple examples show that a direct analogue
of the isotropic formula for the H\"older exponent does not hold in the anisotropic case. Moreover, unless $M$ is isotropic,
this formula never holds for Lipschitz refinable functions, see section \ref{ss.hold.spec}.

Nevertheless, there are ways of treating the anisotropic case. In \cite{Bownik, CGV, CoiW, Schmeisser_Triebel}, the authors consider
special anisotropic H\"older and
Sobolev spaces. We put emphasize on incorporating the spectral properties of the dilation matrix $M$ into the
expression for the H\"older exponent of $\varphi$. Furthermore, we get rid of the $M$-adic expansions
and base our analysis on geometric properties of tilings generated by $M$.

Our paper is organized as follows. In section~\ref{s.hold}, we characterize the continuity and determine the H\"older regularity of
multivariate refinable functions, see Theorems~\ref{th.holder} and~\ref{th.holder-direct}. In subsection \ref{ss.hold.u}, we provide an
algorithm for construction of continuous solutions of \eqref{eq.ref}. We consider several examples and list several important
special cases of Theorems~\ref{th.holder} in subsection~\ref{ss.hold.spec}. The crucial steps of the proofs and actual proofs of
Theorems~\ref{th.holder} and~\ref{th.holder-direct} are given  in subsections~\ref{ss.hold.ideas} and \ref{ss.hold.proofs}.
We
illustrate our results on two numerical examples in subsection~\ref{ss.examples.C}. Both examples deal with continuity
and H\"older regularity of refinable functions, the first one being of shearlet type \cite{CoSauer}. 
In section~\ref{s.higher}, we show how to factorize smooth refinable functions
and compute the H\"older exponents of their directional derivatives.
In section~\ref{s.p}, we analyze the existence of $L_p$-solutions of \eqref{eq.ref}. We show that a direct analogue of
the formula for the H\"older exponent (i.e. replacing the joint spectral radius by the $p$-radius)
does not hold in~$L_p$, $1 \le p < \infty$. To characterize the $L_p$ H\"older exponent of $\varphi$, we consider
extended transition matrices, see subsections \ref{ss.p.Hoelder} and \ref{ss.p.ex}.

\section{Background and notation}\label{s.not}

\noindent We use the standard notation for the function spaces $C, C^k, L_p$, $1 \le p < \infty$.
The space of vector-valued functions $f: X \to \re^n$ with components in $L_p$
is denoted by $L_p(X, \re^n)$. We write $L_p(X)$, if the range space is fixed.
The Schwartz space of smooth rapidly decreasing functions over $\re^s$ is denoted by $\cS$, and $\cS'$ is the space of tempered distributions (distributions over $\cS$ or distributions of
slower growth). By $\mu(X)$ we denote the Lebesgue measure of a set $X \subset \re^n$ and by by $|\cdot |$ either a modulus of a complex number or the cardinality of a finite set. The norm $\|\cdot\|$ in finite dimensional spaces is always а Euclidean, unless stated otherwise.

\subsection{Spectral properties of the dilation matrix} \label{ss.spectral.M}

We assume that the integer dilation matrix~$M \in \z^{s \times s}$ is expansive, i.e.,  all its eigenvalues are
larger than $1$ in the absolute value. Hence, $m=|{\rm det}\, M| \ge 2$. Exactly $n_i$ among the eigenvalues $|\lambda_1| \le \cdots \le |\lambda_s|$ of $M$ are in the absolute value equal to $r_i$, $i=1,  \ldots, q(M) \le s$.
If  $M$ is isotropic, then $q(M)=1$. For $i=1, \ldots , q(M)$, let $J_i\subseteq \re^s$, $\hbox{dim}(J_i)=n_i$, be the eigenspaces of~$M$
corresponding to the eigenvalues of modulus~$r_i$.
The space $\re^s$ is a direct sum
$$
 \re^s\quad =\quad \bigoplus_{i=1}^{q(M)} \, J_i\,
$$
of the  subspaces~$J_1, \ldots , J_{q(M)}$. There exists an invertible transformation $B: \re^s \rightarrow \re^s$ such
that $M$ has the following block diagonal structure
\begin{equation}\label{eq.M-factor}
B^{-1} M B\quad = \quad
\left(
\begin{array}{cccc}
 M|_{J_1} & 0 & \cdots & 0\\
 0       & M|_{J_2} &      & \vdots \\
\vdots   &         & \ddots & 0 \\
0 & \cdots  & 0 & M|_{J_q}
\end{array}
\right),
\end{equation}
where the operator~$M|_{J_i}$ has all its eigenvalues equal to $r_i$ in the absolute value.

\subsection{Dilation matrix and tiles}

The matrix $M$ splits the integer lattice $\z^s$ into $m$ equivalence (quotient) classes defined by
the relation $x \sim y \, \Leftrightarrow \, y - x \in M\, \z^s$.
Choosing one representative $d_i \in \z^s$ from each equivalence class,  we obtain
a {\em set of digits}~$D(M) = \{d_i \ : i=0, \ldots, m-1\}$.  We always assume that~$0 \in D(M)$.
The standard choice is to take $D(M) = \z^s \cap M [0, 1)^s$.

\noindent For every integer point $d \in \z^s$,  we denote by $M_d$, the affine operator
$M_d\, x \, = \, Mx - d, \, x \in \re^s$. We use the notation $0.d_1d_2 \ldots =
\displaystyle \sum_{i=1}^{\infty} M^{-i}d_i$, $d_i \in D(M)$. Consider the following set
\begin{equation}\label{eq.G}
G \quad = \quad \left\{\,\sum_{i=1}^{\infty} M^{-i}d_i  \quad : \quad
d_i \in D(M) \right\}.
\end{equation}
By~\cite{GroH,GroM}, for every expansive
integer matrix~$M$ and for an arbitrary set of digits $D(M)$,
the set $G$ is a compact set with a nonempty interior and possesses  the properties:

\textbf{a)} the  Lebesgue measure $\mu(G) \in \n$;

\textbf{b)}  $\ G \, = \, \bigcup\limits_{d \in D(M)} M_d^{\, -1} \, G$,
the sets $M_d^{\, -1} \, G$ have intersections of zero measure;

\textbf{c)}  the indicator function $\chi = \chi_{G}(x)$ of $G$ satisfies the refinement equation
$$
\chi(x)\  = \  \sum_{d \in D(M)} \, \chi(Mx - d)\, \quad x \in \re^s ;
$$

\textbf{d)}  $\displaystyle \sum_{k \in \z^s} \chi(x+k)\, \equiv \, \mu(G)$, i.e., integer shifts of $\chi$
cover $\re^s$ with $\mu(G)$ layers;
   \smallskip

\textbf{e)}  $\mu(G)=1$ if and only if the function system $\{\chi(\cdot + k)\}_{k \in \z^s}$
is orthonormal.

\noindent If $\mu(G) = 1$, then $G$ is called a {\em tile}.
The integer shifts of a tile define a  {\em tiling}.

\begin{defi}\label{d.tiling}
A tiling generated by  an integer expansive  matrix $M$ and by a set of digits $D(M)$
is a collection of sets $\cG = \{k+ G\}_{k \in \z^s}$ such that

a)  the union of the sets in $\cG$ covers $\re^s$ and $\mu\left( ( \ell+ G) \cap (k+ G)\right)=0$, $\ell \not= k$;

b)  $\displaystyle \ G \, = \, \cup_{d \in D(M)} M_d^{\, -1} \, G$.
\end{defi}
Not every $M$ possesses a digit set $D(M)$ such that $G$ is a tile.  Those situations, however, are rare.
For instance, a digit set generating a tile always exists in cases $s =2, 3$  and also  for arbitrary $s$ with an extra assumption $\, |{\rm det}\,
M|  >  s$, which is quite general for integer expansive matrices~\cite{LW1}. See~\cite{CHM2, LW2} for more details.   Thus, in this paper,
we assume that $G$ is a tile.

\noindent We denote
$$
 G_{d_1 \ldots d_n} \, = \, M_{d_1}^{-1} \cdots M_{d_n}^{-1}G, \quad d_1, \ldots, d_n \in D(M).
$$
Then $\cG^n \, = \, M^{-n}\cG = \{M^{-n}(k+ G)\}_{k \in \z^s}= \, \{G_{d_1 \ldots d_n} \ : \  d_1, \ldots, d_n \in D(M)\}$, $n \in \n$.

\subsection{Refinable functions and the transition operator}

A compactly supported distribution $\varphi \in \cS'(\re^s)$ satisfying equation~(\ref{eq.ref}) is
called {\em a refinable function}. It is well known that
the solution of (\ref{eq.ref}) such that $\int_{\re^s} \varphi (x)\, dx \ne 0$ exists if and only if
$\displaystyle \sum_{k \in \z^s} c_k = m$.   We assume further that the coefficients of (\ref{eq.ref}) satisfy
sum rules of order one
\begin{equation} \label{def:sum_rules}
      \sum_{k \in \z^s} \, c_{Mk-d} \quad = \quad 1\, , \qquad  d \in D(M)\, .
\end{equation}
These conditions are necessary for existence of stable refinable functions~\cite{CDM}.
Consider the {\em transition operator }~$\bT:\, \cS'(\re^s) \, \to \, \cS'(\re^s)$ defined by
\begin{equation}\label{eq.transit}
\bT f \, (x) \ = \ \sum_{k \in \z^s} c_k \, f \, (Mx - k), \quad x \in \re^s.
\end{equation}
For every compactly supported function $f \in \cS'$
such that $\int_{\re^s}f(x)\, dx = 1$, the sequence $\{\bT^jf\}_{j \in \n}$
converges to $\varphi$ in the space~$\cS'$~\cite{CHM2}. The space of distributions
supported on the set
\begin{equation}\label{eq.K}
 K \quad =\quad \Bigl\{\ x \in \re^s \quad : \quad x\ =\ \sum_{j=1}^\infty \
 M^{-j} \gamma_j, \ \  \gamma_j \in \hbox{supp}(c), \ \  c=\{c_k\}_{k \in \z^s}  \Bigl\}
\end{equation}
is invariant under $\bT$. Hence, for $f \in \cS'(K)$, we have $\bT^jf \in  \cS'(K)$ for all~$j \in \n$.
Therefore, the limit~$\varphi \in \cS'(K)$. Thus,  ${\rm supp}\, \varphi \, \subseteq \, K$, see \cite[Proposition 2.2]{CHM2}.

\begin{defi}\label{d.omega}
A finite set $\Omega \subset \z^s$ is a minimal subset of $\z^s$ with the property
$$
 K \quad \subseteq \quad \Omega\ +\ G \ = \ \bigcup\limits_{k \in \Omega}(k + G)\, .
$$
We denote $\, N \, = \, |\Omega|$.
\end{defi}
It is shown easily that $M_d^{-1} (\Omega+G) \, \subseteq \, \Omega+G$ for every $d \in D(M)$.
\smallskip

\noindent The main idea of the matrix approach is to pass from a function~$f: \re^s \to \re$ supported on~$K$
to the vector-valued function
\begin{equation}\label{eq.v}
v:G \rightarrow \re^N, \quad v(x) \ = \ v_{f}(x) \ = \ \Bigl(\, f(x+k)\, \Bigr)_{k \in \Omega}\, , \quad x \in G\, .
\end{equation}
Then the transition operator~(\ref{eq.transit}) restricted to the space
$$
\Bigl\{\, f \in L_1(\re^s) \quad : \quad {\rm supp}\, f\, \subset \, \Omega+G\,  ,
\Bigr\}
$$
becomes the {\em self-similarity operator}~$\bA: L_1(G,\re^N) \rightarrow L_1(G,\re^N)$ defined by
\begin{equation}\label{eq.ss}
(\bA v) (x) \ = \ T_d \, v(Mx - d)\ , \quad x \in M^{-1}(d + G)\ , \quad d \in D(M)\, ,
\end{equation}
where $T_d$ are the $N\times N$ {\em transition matrices} defined by
\begin{equation}\label{eq.T}
(T_d)_{a b} \ = \ c_{\, Ma - b + d}\ , \quad a, b \, \in \, \Omega\, , \quad d \in D(M).
\end{equation}
The rows and columns of the matrices~$T_d$ are enumerated by elements from the set $\Omega$.
We denote
\begin{equation} \label{eq:cT}
 \cT = \{T_d \ : \ d \in D(M)\}.
\end{equation}
The refinement equation becomes the
{\em self-similarity equation} $\bA v = v$ for the vector-valued function~$v(x) = v_{\varphi}(x)$
defined by~(\ref{eq.v}) with~$f = \varphi$, i.e.
\begin{equation}\label{eq.ss1}
v (x) \ = \ T_d \, v(Mx - d)\ , \qquad x \in M^{-1}(d + G)\ , \quad d \in D(M)\, .
\end{equation}

\subsection{Important subspaces of $\re^N$}

We consider the following affine subspace of the space~$\re^N$
$$
 V\ =\ \Bigl\{\, w\, =\, (w_1, \ldots, w_N)^T \in \re^N \quad : \quad \sum_{j=1}^{N} w_j\, =\, 1\, \Bigr\}.
$$
It is well known \cite{CHM2, CDM} that every compactly supported refinable function $\varphi \in \cS'$
such that $\int_{\re^s} \varphi(x)\, ds = 1$ possesses the
{\em partition of unity property:}
$$
\sum_{k \in \z^s} \varphi (x+k)\quad  \equiv \quad  1 \quad x \in \re^s.
$$
Hence, if $\varphi$ is continuous, then $v(x) \in V$ for all $x \in G$.  In particular $v(0) \in V$. For summable refinable function, $v(x) \in V$ for almost all $x \in G$. We denote  the linear part of $V$ by
$$
W\quad =\quad \Bigl\{\, w\ =\ (w_1, \ldots, w_N)^T \in \re^N \quad : \quad  \sum_{j=1}^N w_j=0 \,
\bigr\}\, .
$$
Finally, every continuous refinable function defines the difference space
\begin{equation}\label{eq.U}
 U\ = \ {\rm span}\, \Bigl\{\, v(y)\, - \, v(x) \quad :  \quad  x,y \, \in G\, \Bigr\}, \quad n={\rm dim}\, U.
\end{equation}
Since $v(x) \in V$ for all $x \in G$, we  have  $U \subseteq W$, and, therefore,
$n \le N-1$. The sum rules \eqref{def:sum_rules} imply that the column sums of each matrix $T_d$ are equal to one.
Therefore, $T_dV \subseteq V$ and $T_d W \subseteq W$. Thus, $V$ is a common affine invariant subspace of the family
$\cT$ and $W$ is its common linear invariant subspace.

\noindent Since $U$ is invariant under all $T_d, \, d \in D(M)$, the restrictions $A_d = T_d|_{U}$
of the operators $T_d$, $d \in D(M)$, to the subspace~$U$ are well defined. For a  fixed basis of~$U$, we denote by
\begin{equation} \label{def.cA}
 \cA \ = \ \cT|_U=\{A_d\ :  \ d \in D(M)\}
\end{equation}
the set of the associated $n\times n$ matrices. If the family $\cT$ is irreducible on~$W$, i.e., the
operators in $\cT$ do not share a common nontrivial invariant subspace of W, then
$\cA = \cT|_{W}$.
We also consider the following subspaces of the space~$U$
\begin{equation}\label{eq.ui}
U_i\ = \
{\rm span}\, \Bigl\{v(y)\, - \, v(x) \quad :  \quad  x, y \in G, \  y-x \in J_i \Bigr\}\, , \quad
i=1, \ldots, q(M).
\end{equation}
Note that $U_i$ are nonempty, due to the interior $\hbox{int}(G)$ of $G$ being nonempty. It is seen easily that the spaces $\{U_i\}_{i=1}^{q(M)}$ span the whole space $U$,
but their sum may not be direct. Indeed, the subspaces $\{U_i\}_{i=1}^{q(M)}$, unlike the subspaces $\{J_i\}_{i=1}^{q(M)}$, may have nontrivial intersections. For example, they can all coincide with $U$. Lemma \ref{l.invar} shows that all $U_i$ are invariant under $\cA$.

\begin{lemma}\label{l.invar}
If $J$ is an invariant subspace of $M$, then
$L \, = \, {\rm span}\, \{v(y) - v(x) \ : \ y-x \in J\}$ is a common invariant subspace for $\cA$.
\end{lemma}
\begin{proof} If $u \in L$, then $u$ is a linear combination
of several vectors of the form $v(y) - v(x)$ with $y - x  \in J$.
For $x' = M^{-1}(x+d)$, $y' = M^{-1}(y+d)$, $d \in D(M)$, we have $y'-x' \in J$ and
$$
v(y') - v(x') \ = \ A_d \, \bigl( v(My'-d) \, - \, v(Mx' - d)\bigr) \ = \ A_d \, \bigl( v(y) \, - \, v(x)\bigr)\, .
$$
Hence, $A_d\bigl( v(y) -  v(x)\bigr) \in L$ for each pair $(x, y)$, and, therefore, $A_d u \in L$ for all $u \in L$.
\end{proof}

\subsection{Joint spectral radius}

\begin{defi}\label{d.jsr}
The joint spectral radius of a finite family~$\cA$ of linear operators $A_d$ is
$$
\rho(\cA)\ = \ \lim_{k \to \infty}\ \max_{A_{d_i} \in \cA, \, i = 1, \ldots , k}\
\|A_{d_1}\cdots A_{d_k}\|^{\, 1/k}.
$$
\end{defi}
\noindent This limit always exists and does not depend on the operator norm~\cite{RS}.
The joint spectral radius measures the simultaneous contractibility of the
operators from~$\cA$. Indeed, $\rho(\cA) < 1$ if and only if there exists a norm in~$\re^n$
in which all $A \in \cA$ are contractions. In general,
$$
 \rho(\cA)\quad =\quad \inf\, \Bigl\{\, \beta  \ge 0 \ :  \ \exists \ \| \cdot \| \ \hbox{in} \ \re^n \ \hbox{such that} \  \|A\| < \beta, \ A \in \cA \, \Bigr\}.
$$
We denote
$$
\rho_i\  = \ \rho(\cA|_{U_i}), \quad i=1,\ldots,q(M).
$$

\section{Continuous solutions and H\"older regularity}\label{s.hold}

In this section, in Theorem \ref{th.holder}, we characterize the continuity of a solution $\varphi$ of the refinement equation
\eqref{eq.ref}
in terms of the spectral  properties of $\cA$ and determine the exact H\"older exponent
$$
 \alpha_\varphi\ =\ \sup \, \bigl\{\alpha \ge 0 \ : \  \|\varphi(\cdot+h)-\varphi \|_{C(\re^s)} \le C \|h\|^\alpha, \  \ h \in \re^s \bigr\}
$$ of $\varphi$.
Although the definition of $U$ in \eqref{eq.U} depends on
$\varphi$,  Propositions \ref{p.U0} and \ref{p.d.U} remove this dependency. Moreover, the space $U$ can be
determined explicitly using Algorithm~1 in subsection
\ref{ss.hold.u} without the knowledge of $\varphi$. If this
algorithm fails, then there exists no continuous solution of the
corresponding refinement equation. The special cases of Theorem
\ref{th.holder} are considered in subsection \ref{ss.hold.spec},
for its summary see Remark \ref{r.410}. The crucial result for the
proof of Theorem \ref{th.holder} is Theorem
\ref{th.holder-direct}. The main steps of the proof of
Theorem \ref{th.holder-direct} are summarized
in subsection \ref{ss.hold.ideas} and the proofs of Theorems \ref{th.holder} and
\ref{th.holder-direct} are given in subsection \ref{ss.hold.proofs}. We illustrate our results with
examples in subsection \ref{ss.examples.C}. For the readers convenience, we first list some of the above
mentioned crucial results.

\begin{prop}\label{p.U0}
Let $v_0 \in V$ be an eigenvector of $T_0$ associated to the eigenvalue $1$.
If $\varphi \in C(\re^s)$, then $U$ is the smallest by inclusion common invariant
subspace of the matrices $T_d, \, d \in D(M)$, that contains the~$m$ vectors $\, T_dv_0 - v_0, \, d \in D(M)$.
\end{prop}

\noindent For the proof of Proposition \ref{p.U0} see the arXive version of this paper \cite{CharinaProtasov}. Our proof is the
multivariate analogue of \cite[ Proposition 3]{CH}.

\begin{remark}\label{r.300} {\em Recall that $0  \in D(M)$, which justifies the notation~$T_0$.
The existence of the eigenvector $v_0 \in V$ of $T_0$ associated to the eigenvalue~$1$
follows by by~(\ref{eq.ss1}), i.e. $T_0v(0)\, = \, v(0)$ and by the continuity $\varphi$ implying that
$v(0) \in V$.}
\end{remark}

\noindent Similarly to  \cite{CH}, Proposition~\ref{p.U0} yields a useful characterization of $U$.

\begin{prop}\label{p.d.U}
Let $v_0\in  V$ be an eigenvector of $T_0$ associated to the eigenvalue~$1$.
 The space $U$ is the minimal common invariant subspace of the $m$ matrices $T_d, \, d \in D(M)$,
 and $U$ contains the $m$ vectors $T_dv_0 - v_0, \, d \in D(M)$.
\end{prop}


\noindent 
If the eigenvalue $1$ of $T_0$ is not simple,
Proposition~\ref{p.U3} in subsection~\ref{ss.hold.u} guarantees that there exists at most one eigenvector $v_0 \in V$ associated with
one of these eigenvalues $1$ such that
$\rho(\cA)<1$.  Thus, the subspace $U$ can be computed, unless the refinement equation does not possess a continuous solution.
For simplicity, we make the following assumption.

\begin{assum}\label{a.1}
The matrix $T_0$ has a simple eigenvalue $1$.
\end{assum}

\noindent Recall that $\rho_i\  = \ \rho(\cA|_{U_i})$, where the subspaces $U_1, \ldots , U_{q(M)}$ are defined in~(\ref{eq.ui}).

\begin{theorem}\label{th.holder}
A refinable function $\varphi$ belongs to $C(\re^s)$ if and only if $\rho(\cA) < 1$. In this case,
\begin{equation}\label{eq.holder}
\alpha_{\,\varphi}\quad = \quad   \min\limits_{i = 1, \ldots , q(M)}\, \log_{\, 1/r_i} \, \rho_i\, .
\end{equation}
\end{theorem}

\noindent The proof of~(\ref{eq.holder}) is based on Theorem~\ref{th.holder-direct}. To state it, we define
the {\em H\"older exponent of $\varphi$ along a linear subspace~$J \subset \re^s$} by
$$
 \alpha_{ \varphi , J} \quad  =\quad  \sup \, \bigl\{\alpha \ge 0 \ : \ \|\varphi(y) - \varphi(x)\| \le C \|y - x\|^{\, \alpha}\, , \
y-x \in J\, \bigr\}\, .
$$

\begin{theorem}\label{th.holder-direct}
If $\varphi \in C(\re^s)$, then
\begin{equation}\label{eq.holder-direct}
\alpha_{\, \varphi ,  J_i}\quad = \quad  \log_{1/r_i} \, \rho_i, \quad i = 1, \ldots , q(M)\, .
\end{equation}
\end{theorem}

\begin{remark}\label{r.410}
{\em The identity \eqref{eq.holder} emphasizes the influence  of the spectral structure of the
dilation matrix $M$ on the regularity of the solution~$\varphi$.
Recall that, in the univariate case, the H\"older exponent is given by  $\alpha_{\varphi} \, = \, \log_{1/r}\, \rho(\cA)$, where $M=r \ge 2$ is the corresponding dilation factor.
In the multivariate case, the H\"older exponent is equal to the minimum of several
such values taken over different dilation coefficients~$r_i$ on the corresponding subspaces~$J_i$
of~$M$. In special, favorable multivariate cases, the expression in \eqref{eq.holder} becomes $\alpha_\varphi=\log_{1/\rho(M)}\, \rho(\cA)$
and, thus, resembles the univariate case. This happens, for instance,  when the matrix $M$ is isotropic, i.e. $|\lambda_1|= \ldots= | \lambda_s| = \rho(M)$, in particular, when
 $M = r\, I$, $r \ge 2$. Another favorable situation is when the matrices in $\cA$
do not possess any common invariant subspace. However, the need for the minimum in~(\ref{eq.holder}) is not exceptional. It is of
crucial importance e.g. for anisotropic refinable Lipschitz continuous functions $\varphi$,
see Corollary~\ref{c.regular} in subsection~\ref{ss.hold.spec}.}
\end{remark}

\subsection{Special cases of Theorem \ref{th.holder} and examples}\label{ss.hold.spec}

\noindent To compare the result of Theorem \ref{th.holder} with the known results from the
literature, we need to define the stability of $\varphi$.

\begin{defi}\label{d.stable}
A compactly supported $f \in L_\infty(\re^s)$ is stable, if there
exists $0<C_1 \le C_2 <\infty$ such that for all $\bc \in
\ell_\infty(\z^s)$,
$$
 C_1 \|\bc\|_{\ell_\infty} \le \big\| \sum_{\alpha \in \z^s} \bc(\alpha)f(\cdot-\alpha)
 \big\|_\infty \le C_2\|\bc\|_{\ell_\infty}\,.
$$
\end{defi}

\noindent \textbf{The univariate case} ($\mathbf{s=1}$). In this case, the dilation factor is $M\ge 2$ and $M=m=r$. Theorem~\ref{th.holder} becomes a well-known statement that
$\alpha_{\varphi} = \log_{1/r} \rho(\cA)$. If $\varphi$ is stable,  then we have~$\rho(\cA)=\rho(\cT|_{U}) = \rho(\cT|_{W})$ even if $U \ne W$ (see~\cite{CDM}). The space~$U$ was completely characterized in~\cite{P07} and it was shown that
every refinement equation can be factorized to the case $U = W$. In the multivariate case, however, there is no general factorization procedure and,
even in the stable case, we cannot expect~$U = W$, see Example~\ref{ex.10} below.
\bigskip

\noindent \textbf{The case  $\mathbf{s \ge 2}$ with isotropic dilation matrix}. Since $q(M) = 1$, it follows that $U_1 = U$. Theorem~\ref{th.holder} then implies the following well-known fact.
\begin{cor}\label{c.iso}
If $M$ is isotropic, then $\alpha_{\varphi}\, = \, \log_{\, 1/\rho(M)} \rho(\cA)$.
\end{cor}
\bigskip

\noindent \textbf{The irreducible case with} $\mathbf{s \ge 2}$. The dilation matrix $M$ can be anisotropic, i.e. the number of
different in modulus eigenvalues of $M$ is $q(M)>1$.
We say that the set of matrices~$\cA = \cT|_{\, U}$ is irreducible, if they do not possess any common invariant subspace.
Another corollary of Theorem~\ref{th.holder} states the following.
\begin{cor}\label{c.irred}
If the family $\cA$ is irreducible, then $\alpha_{\varphi}\, = \,  \log_{\, 1/\rho(M)} \rho(\cA)$.
\end{cor}

\noindent The irreducibility assumption fails however in many important cases. For instance, if~$\varphi$
is a tensor product of two refinable functions, then $\cA$ is always reducible.

\begin{ex}\label{ex.10}
{\em Let $\varphi_1, \varphi_2 \in C^1(\re)$ be two univariate refinable function with dilations $M_1=2$ and $M_2=3$ and
refinement coefficients $\bc_1 \in \ell_0(\z)$ and $\bc_2 \in \ell_0(\z)$, respectively. Then the function
$\varphi = \varphi_1 \otimes \varphi_2$ satisfies the refinement equation
with $M = {\rm diag}\, (2\, , \, 3)$ and ~$\bc = \bc_1\otimes\bc_2$. Due to $\varphi_1, \varphi_2 \in C^1(\re)$, we have
$\rho_1 = \rho(\cA|_{U_1}) = \frac12, \ \rho_2 = \rho(\cA|_{U_2}) = \frac13$.
By Theorem~\ref{th.holder}, $\alpha_{\varphi} \, = \, \min\, \{ \log_{1/2} \rho_1\, , \, \log_{1/3} \rho_2\} = 1$,
which is natural, because $\varphi \in C^1(\re^2)$. On the other hand,
$\rho (\cA) = \max \, \bigl\{ \rho_1\, , \,
\rho_2 \bigr\} = \frac12$. Hence, $ \log_{\, 1/\rho(M)} \rho(\cA) =  \log_{1/3} \frac12 \, = \,
0.630092\ldots $.
Thus, in this case, $\alpha_{\varphi}\, >  \,  \log_{\, 1/\rho(M)} \rho(\cA)$.
Note that, if $\varphi_1$ and $\varphi_2$ are both stable, then so is $\varphi$.
Nevertheless, unlike in the univariate case, the H\"older exponent of $\varphi$ is not determined by
the value $\log_{\, 1/\rho(M)} \rho(\cA)$.
 }
\end{ex}

\noindent After  Example~\ref{ex.10} one may hope that
the case of reducible family~$\cA$ is exceptional, and the equality
$\alpha_{\varphi} = \log_{\, 1/\rho(M)} \rho(\cA)$ actually holds for most  refinable functions.
On the contrary, the  result of Corollary~\ref{c.regular} shows that the the situation when the isotropic  formula fails is rather generic.

\begin{cor}\label{c.regular}
If the matrix $M$ is anisotropic and the refinable function $\varphi \not =0$
is Lipschitz continuous, then $1=\alpha_\varphi> \log_{1/\rho(M)} \rho(\cA)$ and the family~$\cA$ is reducible.
\end{cor}
\begin{proof} Assume that $1\, \le  \,  \log_{\, 1/\rho(M)} \rho(\cA)$,  or, equivalently,  $\rho(\cA) \le  1/\rho(M)$. Since $M$ is anisotropic, factorization~(\ref{eq.M-factor})
contains $q(M)\ge 2$ blocks, and, hence, $r_i < \rho(M)$ for some $i \in \{1, \ldots, q(M)\}$. By Theorem~\ref{th.holder-direct},
we have $\alpha_{\varphi, J_i}\, = \, \log_{1/r_i} \rho(\cA) \, > \, \log_{1/\rho(M)} \rho (\cA) \, \ge \, 1$. Therefore,
$\varphi$ is constant on every affine subspace $u + U_i$, $u \in \re^n$.
Hence, $\varphi \equiv 0$, because it is compactly supported. The reducibility of $\cA$ follows by Corollary \ref{c.irred}.
\end{proof}

\noindent Thus, we see that at least for all anisotropic smooth refinable functions, the simple formula
for the H\"older exponent fails and the minimum in~(\ref{eq.holder}) is significant.

\noindent \textbf{The case of a dominant invariant subspace}. In practice, this case is much more generic than the irreducible case.
\smallskip

\begin{defi} \label{d.dominant} A subspace $U' \subset U$ is called {\rm dominant} for a family of operators $\cA$ if
\begin{description}
\item[$(i)$] $U'$ is a common invariant subspace of $\cA$,
\item[$(ii)$] $U'$ is contained in all common invariant nontrivial subspaces of~$\cA$ and
\item[$(iii)$] $\rho(\cA|_{U'})\, = \, \rho(\cA)$.
\end{description}
\end{defi}

\noindent Take a basis of a dominant subspace $U'$ and complement it to a basis  of $U$. Let $B$ be the $n \times n$ matrix containing
these basis elements of $U$.
Then every matrix $A_d \in \cA$ in this basis has a block lower triangular form
\begin{equation}\label{eq.domin}
B^{-1}A_d\,B \quad = \ \quad
\left(
\begin{array}{cc}
\tilde{A_d} & 0\\ * & A_d|_{U'} \end{array} \right), \quad d \in D(M).
\end{equation}
By Definition \ref{d.dominant},
$$
\rho(\cA|_{U'})\  = \ \max\, \bigl\{\, \rho(\tilde{\cA}), \rho(\cA|_{U'})\bigr\}\ = \ \rho(\cA)\, , \quad
\tilde{\cA}=\{\tilde{A}_d\ : \  d \in D(M)\}.
$$
Furthermore, since any common invariant subspace of $\cA$ contains $U'$, it follows that
the joint spectral radius of $\cA$ restricted to any common invariant subspace is equal to
$\rho(\cA)$. Therefore, we have proved the following result.
\begin{cor}\label{c.domin}
If the family $\cA$ possesses a dominant subspace, then $\, \alpha_{\varphi}\, = \,  \log_{\, 1/\rho(M)} \rho(\cA)$.
\end{cor}

\subsection{Construction of $U$ and of the continuous refinable function~$\varphi$.}\label{ss.hold.u}

\noindent In this section, we answer two crucial questions: how to determine the space $U$
and how to construct  the corresponding continuous refinable function~$\varphi$.
In the univariate case, the algorithm for determining the space~$U$ was elaborated in~\cite{CH}.
In this section, we present its multivariate analogue.

\medskip
\noindent \textbf{Algorithm for construction of the space~$U$.}

\vspace{0.3cm} \noindent {\bf Algorithm 1:} For a given set $\cT=\{T_d \in \re^{N \times N}\ : \ d \in D(M)\}$ of transition matrices

\begin{description}

 \item{\it $1$.Step:} Compute $v_0$ in  $T_0 \, v_0=v_0$ and normalize $( {\bf 1}, v_{0})\, =\, 1$, where ${\bf 1} = (1, \ldots , 1)^T \in \re^N$.

 \item{\it $2$.Step:} Define $U^{(1)}={\rm span} \, \bigl\{T_dv_0\, -\, v_0 \ : \ d \in D(M) \setminus \{0\}\bigr\}$.

 \item{\it $3$.Step:} Repeat
    $$
          U^{(k+1)}= {\rm span} \, \bigl\{ U^{(k)},  T_d u^{(k)} \ : \ u^{(k)} \in U^{(k)}, \quad d \in D(M)\},
            \quad 1 \le k \le N-1,
    $$
    {\bf while} ${\rm dim}(U^{(k)}) \ <  \ {\rm dim}(U^{(k+1)})$.
 \end{description}

\noindent \begin{equation} \label{algor:def:U}
\hspace{-12cm} {\rm Output:} \quad  U=U^{(k)}
\end{equation}

\noindent If {\it $1$.Step} is impossible, i.e., the eigenvector~$v_0$ does not exist, then, by Remark~\ref{r.300},
the solution of the refinement equation \eqref{eq.ref} is not continuous. Note that $1 \le k \le N-1$ is dictated by
${\rm dim}\, U \, \le \,  N-1$ and that, by construction, at least one extra element is added to $U^{(k+1)}$ before the algorithm
terminates. Proposition~\ref{p.U0}, stated at the beginning of section \ref{s.hold}, implies that the space $U$ in Algorithm 1 coincides
with the space $U$ in Definition~\ref{eq.U}.

\medskip
\noindent \textbf{Algorithm for construction of a continuous $\varphi$}.
\smallskip

\noindent Due to the fact that the rational $M$-adic points are dense in $G$, the slight modification of Algorithm~1
yields a method for the step-by-step construction of the vector-valued function $v = v_{\varphi}$ defined on $G$ or,
equivalently, of the function~$\varphi$.

\vspace{0.3cm} \noindent {\bf Algorithm 2:} For a given set $\cT=\{T_d \in \re^{N \times N}\ : \ d \in D(M)\}$ of transition matrices

 \begin{description}
 \item{\it $1$.Step:} Compute $(0,v_0)$ such that  $T_0 \, v_0=v_0$ and normalize
  $( {\bf 1}, v_0)\, = \, 1$.

 \item{\it $2$.Step:} Define $V^{(0)}\, =\, \{ (0, v_0)\}$.

 \item{\it $3$.Step:} For $k=1,2,3, \ldots$
    $$
          V^{(k)}= V^{(k-1)} \ \cup \ \bigl\{ (x,v_x) \ : \
          x\, =\, M_d^{-1} y, \  v_x\, =\, T_d v_y, \ (y,v_y) \in V^{(k-1)},
            \  d \in D(M)\}
        $$
        \hspace{0.5cm}end
 \end{description}

\noindent If the function $v$ is continuous, Algorithm~2 determines $V^{(k)}$
in a unique way.

\noindent The next result ensures that $U$ is well defined even if the eigenvalue $1$ of the matrix~$T_0$ is not simple.

\begin{prop}\label{p.U3} The matrix $T_0$ has at most one eigenvector $v_0 \in V$
associated with the eigenvalue $1$ such that $\rho(\cA) < 1$.
If such $v_0 \in V$ exists, then
$\varphi$ is continuous and $v_0 = v_{\varphi}(0)$.
\end{prop}
\begin{proof} By Algorithm 2, there exists a refinable function $\varphi$ with $v_0 = v_\varphi(0)$.
If there is another eigenvector $\tilde{v}_0 \in V$ with this property, then, by Algorithm 2, it generates another refinable function
$\tilde{\varphi}$ for which $\tilde{v}_0 = v_{\tilde{\varphi}}(0)$. By the uniqueness of the solution of the refinement equation,
these two solutions may only differ by a constant, hence, the vectors $v_0$ and $\tilde{v}_0$ are collinear.
The continuity of $\varphi$ follows from Theorem~\ref{th.holder}.
\end{proof}

\subsection{Road map of our main results}\label{ss.hold.ideas}

\noindent We would like to emphasize that, to tackle the anisotropic case,
we use geometric properties of tilings rather than the $M$-adic
expansions of points in~$\re^s$ (the latter being a successful
strategy in the isotropic case). Our key contribution is
Theorem~\ref{th.holder-direct} that finally reveals the delicate
dependency of the H\"older exponent of a refinable function on its
H\"older exponents along the subspaces $J_i$, $i=1, \ldots,q(M)$. Due to the
importance of Theorem~\ref{th.holder-direct}, we would like to
give here a preview of its proof.
\smallskip

{\em Step 1.} Extend the vector-valued function~$v$ in
\eqref{eq.v} defined on the tile $G$ to the whole~$\re^s$,
see~(\ref{d.tildev}). Lemma~\ref{l.main} yields, for $x, y \in
\cG$ (i.e. $x,y \in G-j$ for some $j \in \n$), the estimate
$\|\tilde v(y) - \tilde v(x)\|\, \le \, \| v(y+j) - v(x+j) \|$.
The extension of $v$ is motivated by the fact that parts of the
line segment $[x,y]$ can lie outside of $G$, due to its possible
fractal structure.
\smallskip

{\em Step 2.} Lemma~\ref{l.tiling} shows that, for a tiling $\cG$ of~$\re^s$,
the total number of the subsets of the tiling intersected by a line segment is proportional to the length of that segment.
\smallskip

{\em Step 3.} Due to Step 2, Lemma~\ref{l.colours} and Proposition~\ref{p.aux} imply
that, for $k \in \n$, any line segment $[x, y]$ in $\re^s$, $y - x \, \in \, J_i$, consists of several line segments
such that 1) the endpoints of each of those line segments belong to one subset of the tiling~$\cG^k$; 2) the total number of
those line segments is bounded by $C\, \|M^k(x - y)\|\, \asymp \,  r_i^k \|x-y\|$.
\smallskip

{\em Step 4}. The difference between the values of the function $v$ at the endpoints of each
of those subsegments  of $[x,y]$ is bounded from above by $C_1\, (\rho_i + \varepsilon)^k$ for some $\varepsilon>0$.
Hence, by Step 1, the same is true for~$\tilde v$.
Therefore, by the triangle inequality,
$\|\tilde v(y) - \tilde v(x)\|\, \le \, C_2 (\rho_i + \varepsilon)^k$.
For $k$ such that $r_i^k \, \asymp \, 1 / \|x-y\| $, we obtain
$\|\tilde v(y) - \tilde v(x)\|\, \le \, C_3 \|x-y\|^{\, \alpha (\varepsilon) }$,
where $\alpha(\varepsilon)$ approaches $\log_{1/r_i} \rho_i$ as $\varepsilon$ goes to $0$.

\subsection{Auxiliary results for Theorems~\ref{th.holder} and \ref{th.holder-direct} }\label{ss.hold.aux}

\noindent The proofs of our main results, Theorems~\ref{th.holder} and \ref{th.holder-direct}, are based on an important observation
formulated in Proposition~\ref{p.aux}. We also make use of the following basic properties of the joint spectral radius
and two auxiliary lemmas.

\smallskip

\noindent \textbf{Theorem A1}~\cite{RS}. {\it For a family of operators~$\cA$ acting in~$\re^n$ and for any
$\varepsilon > 0$, there exists a norm $\|\cdot \|_\varepsilon$ in $\re^n$  such that
$\|A\|_\varepsilon \, < \, \rho(\cA) + \varepsilon$ for all $A \in \cA$.}

\medskip

\noindent \textbf{Theorem A2}~\cite{B}. {\it For a family of operators~$\cA$ acting in~$\re^n$
there exists $u \in \re^n$ and a constant
$C(u) > 0$ such that $\max_{A_{d_i} \in \cA}\|A_{d_1}\cdots A_{d_k}u\| \ge C(u)\, \rho(\cA)^{\, k}, \ k \in \n$.
Moreover, if $\cA$ is irreducible, then $\max_{A_{d_i} \in \cA}\|A_{d_1}\cdots A_{d_k}\| \le  C\, \rho(\cA)^{\, k}, \  k \in \n$,
for some constant $C>0$. }

\begin{lemma}\label{l.colours}
Assume that the segment $[0, 1]$ is covered with $\ell$ distinct closed sets.
Then there exist $\ell+1$ points $0 = a_0 \le \ldots \le a_{\ell} = 1,$
such that for each $i = 0, \ldots , \ell-1$,  the points $a_i, a_{i+1}$ belong to one of these sets,.
\end{lemma}
\begin{proof} Let the first set contain the point $a_0=0$.
Choose $a_1$ to be the maximal (in the natural ordering of the real line) point
of the first set. If $a_1 \ne 1$, then $a_1$ must belong to another set of the tiling.
Choose $a_2$ to be the maximal point of this set. Repeat until $a_{\ell_0} = 1$
for some~$\ell_0 \in \n$. We have $\ell_0 \le \ell$, since  the sets are distinct.
If $\ell_0 < \ell$, we extend the sequence $a_0 \le \ldots \le a_{\ell_0}$ by  the points $a_{\ell_0+1}, \ldots ,  a_{\ell} = 1$.
\end{proof}

\noindent Next we show that a segment of a given length intersects finitely many sets of the tiling $\cG$.

\begin{lemma}\label{l.tiling}
For a tiling $\cG$, there exists a constant $C>0$ such that
every line segment $[x, y] \in \re^s$ intersects at most $\, C  \max\{1\, , \, \|y-x\|\, \}$
sets of $\cG$.
\end{lemma}
\begin{proof} It suffices to prove that the number of sets $G+k \subset \cG$ intersected by a segment of length $1$ is bounded
above by $C>0$. It will imply that the number of sets $G+k \subset \cG$ intersected by a segment of length $\|y-x\|>1$,
 is bounded by $C\, \|y-x\|$, and the claim follows. Thus, let a segment $[x, y]$ be of length $1$.
If $(G+k) \cap [x,y] \not= \emptyset$ for some $k \in \n$, then the set $G+k$ is contained in
$[x,y]+B_r(0)$, where $B_r(0)$ is the Euclidean ball of radius $r={\rm diam}\, (G)$. Denote by $V$ the volume of $[x,y]+B_r(0)$,
then the total number of sets $G+k \subset \cG$ intersecting $[x, y]$ is bounded by
$C=\frac{V}{\mu(G)} = V < \infty$, due to $\mu(G)=1$.
\end{proof}

\noindent To deal with line segments $[x,y]$, $x,y \in G$, that do not completely belong to $G$, we
extend the continuous vector-valued function $v=v_\varphi$ in \eqref{eq.v} which is defined on $G$ to the whole~$\re^s$.
Define
\begin{equation}\label{d.tildev}
 \tilde{v}: \re^s \rightarrow \re^N, \qquad \tilde{v}(x)=\left( \varphi(x+k) \right)_{k \in \Omega}.
\end{equation}
In Lemma \ref{l.main} and in Proposition \ref{p.aux}, we compare the properties of $v$ and $\tilde{v}$.

\begin{lemma}\label{l.main}
Let $x, y \in  G-j$, $j \in \z^s$. Then
$\|\tilde v(y) \, - \, \tilde v(x)\| \, \le \, \|v(y+j) \, - \, v(x+j)\|$.
\end{lemma}
\begin{proof} Let $j \in \z^s$. By \eqref{eq.v} and due to the compact support of $\varphi$, the $k$-th component
of $\tilde v(y) \, - \, \tilde v(x)$ is given by
$$
  \left( \tilde v(y) \, - \, \tilde v(x)\right)_k=\left\{ \begin{array}{ll}
        \left( v(y+k) \, - \, v(x+k)\right)_\ell, & \ell=-j+k \in \Omega, \\ 0, & \hbox{otherwise},
    \end{array} \right. \quad k \in \Omega.
$$
Hence, we have $\|\tilde v(y)  - \tilde v(x)\| \le \|v(y+j) - v(x+j)\|$.
\end{proof}

\begin{prop}\label{p.aux}
Let $\varphi \in C(\re^s)$ be refinable, $x, y \in \re^s$ and $k \in \n$.
There exist
$$\ell \, \le \, \max\, C \, \bigl\{\, 1\, , \, \|M^k(x-y)\|\, \bigr\}$$ (with $C>0$ from Lemma~\ref{l.tiling}),  integers $\{d_1^{(i)},
\ldots, d_k^{(i)}\}_{i=0}^{\ell-1}$ from $D(M)$,
positive numbers $\{\alpha_i\}_{i=0}^{\ell - 1}$ whose sum is equal to one,
and sets of points $\{x_i\}_{i=0}^{\ell - 1}$,  $\{y_i\}_{i=0}^{\ell - 1}$ from $G$ such that
$y_i - x_i \, = \, \alpha_i\, M^k(y-x)$ for all $i = 0, \ldots , \ell-1$, and
\begin{equation}\label{eq.auxmain}
 \bigl\| \tilde v(y) \, - \, \tilde v(x)  \bigr\|\quad \le \quad
 \sum_{i=0}^{\ell-1} \  \Bigl\| T_{d_1^{(i)}} \cdots T_{d_k^{(i)}}\, \Bigl( v\, \bigl( y_i \bigr)\ -  \ v\, \bigl( x_i \bigr) \Bigl)\, \Bigr\|\, .
\end{equation}
\end{prop}
\begin{proof} For $[x,y] \subset \cG^k$, by Lemma~\ref{l.colours},
there exist $\ell+1$ points $\{x=a_0 \le a_1 \le  \ldots \le a_\ell=y\} \subset [x,y]$ such that each pair of
successive  points~$a_i, a_{i+1}$ belongs to only one set $G_{d_1^{(i)} \ldots d_k^{(i)}}-j^{(i)}$, $j^{(i)} \in \n$,
of the tiling~$\cG^k$. First we give an estimate for~$\ell$. Since $\ell$ elements of the tiling
$\cG^k \, = \, M^{-k}\cG$ cover a segment of length~$\|y - x\|$,
the same number of elements of the tiling $\cG$ cover a segment of length $\|M^k(y-x)\|$.
Therefore,  Lemma~\ref{l.colours} yields   $\, \ell \, \le \, C \, \max\{1\, , \, \|M^k(y-x)\|\, \}$.
Furthermore, the set  $G_{d_1^{(i)} \ldots d_k^{(i)}}-j^{(i)} \subset G-j^{(i)}$, $i=0, \ldots, \ell-1$.
By Lemma~\ref{l.main}, we obtain
$$
\Bigl\|\tilde v(y) \, - \tilde v(x)\Bigr\|\ \le \ \sum_{i=0}^{\ell - 1} \ \Bigl\|\tilde v(a_{i+1}) \, - \, \tilde v(a_i)\Bigr\|\
\le \ \sum_{i=0}^{\ell - 1} \ \Bigl\|\, v\, \bigl(a_{i+1} + j^{(i)}\bigr) \, - \,  v\bigl(a_i + j^{(i)}\bigr)\, \Bigr\|\, .
$$
Due to $a_i + j^{i}\, , a_{i+1}+j^{(i)} \in G_{d_1^{(i)} \ldots d_k^{(i)}}$, $i=0,\ldots,\ell-1$,
the points
$$
x_i \quad = \quad M_{d_k^{(i)}} \cdots M_{d_1^{(i)}}\, (a_i+j^{(i)}) \quad \hbox{and} \quad
y_j \quad = \quad M_{d_k^{(i)}} \cdots M_{d_1^{(i)}}\, (a_{i+1}+j^{(i)})
$$
belong to $G$. Thus, by~(\ref{eq.ss1}), we obtain
\begin{equation}\label{eq.auxmain1}
 \Bigl\|\tilde v(y) \, - \tilde v(x)\Bigr\| \ \le \sum_{i=0}^{\ell-1} \Bigl\| T_{d_1^{(i)}} \cdots T_{d_k^{(i)}}
 \Bigl( v\bigl(y_i) \, \bigr)\  - \  v\bigl(x_i\, \bigr)\,  \Bigr)\, \Bigr\|\, .
\end{equation}
For each $i=0, \ldots, \ell-1$, we define the number $\alpha_i$ from the equality
$\|a_{i+1} - a_i\|\, = \, \alpha_i\, \|y-x\|$.
It follows that $\, \displaystyle  \sum_{i=1}^{\ell - 1}\, \alpha_i \, = \, 1$ and  that
$y_i- x_i \, = \, M^k (a_{i+1} - a_i)\, = \, \alpha_i M^k (y-x)$.
\end{proof}

\subsection{Proofs of Theorems~\ref{th.holder} and~\ref{th.holder-direct}}\label{ss.hold.proofs}

In this subsection we prove Theorems~\ref{th.holder} and \ref{th.holder-direct}.
We start with Theorem~\ref{th.holder-direct} as its proof is a crucial part of the
proof of Theorem~\ref{th.holder}. Note that for both Theorems~\ref{th.holder} and \ref{th.holder-direct}
the assumption that  $\varphi \in C(\re^s)$ implies, e.g. by \cite{CHM2}, that $\rho(\cA)<1$. We will not
reprove this result here.

\begin{proof}[{\tt Proof of Theorem~\ref{th.holder-direct}}] Let $\varepsilon \in (0, 1-\rho_i)$ and
$i \in \{1, \ldots , q(M)\}$. We first show that $\alpha_{\varphi, J_i}
\ge \log_{1/r_i} \, \rho_i$. For arbitrary points $x, y \in G$ such that $y - x \in J_i$ and $\|y - x\| < 1$, define $k$
to be the smallest integer such that $ \, \|M^k(y - x)\| \, \ge \, 1$.  Since $y - x \in J_i$, it follows that
\begin{equation}\label{eq.k0}
1 \le \|M^{k} (y - x)\| \quad \le \quad C\, (r_i+\varepsilon)^k\, \|y - x\|\, ,
\end{equation}
where the constant $C>0$ depends only on~$M$. By
\eqref{d.tildev}, Theorem~A1 and by Proposition~\ref{p.aux}, for
these $x, y$ and $k$, there exist a constant $C_1>0$ depending on
$G$ and the integer
$$
\ell \  \le \  C_1\, \max\, \bigl\{1\, , \, \|M^k(y-x)\|\, \bigr\}\ = \ C_1 \|M^k(y - x)\|
$$
such that (note that $y-x \in J_i$ implies, by
Proposition~\ref{p.aux}, that $x_j,y_j$ in \eqref{eq.auxmain1}
satisfy $y_j-x_j \in J_i$, $j=0,\ldots,\ell-1$)
$$
\bigl\|v(y) \, - \, v(x)\bigr\| \le  2\, \ell\,  C_2 \, (\rho_i+\varepsilon)^k  \, \|v \|_{C(G)} \le
 2\, C_3 \|M^k(y - x)\| \, \|v\|_{C(G)}\, (\rho_i+\varepsilon)^k
$$
with the constant $C_3$ independent of $k$.
By the choice of $k$, we have $\|M^{k-1}(y - x)\| < 1$ and, hence,
$\|M^k(y - x)\| \, \le  \, \|M\|\, \|M^{k-1}(y - x)\|\, < \, \|M\|$. Thus,
\begin{equation} \label{eq.some}
\bigl\|v(y) \, - \, v(x)\bigr\| \ \le \ 2\, C_3
 \, \|M\| \, \|v\|_{C(G)}\, (\rho_i+\varepsilon)^k.
\end{equation}
Combining the above estimate with \eqref{eq.k0} (i.e. $k \ge
-\frac{\log\|y-x\|}{\log(r_i+\varepsilon)}+C_4$), we get, due to
$\rho_i+\varepsilon<1$,
$$
\bigl\|v(y) \, - \, v(x)\bigr\| \ \le \ C(\varepsilon)\, \|y - x\|^{\, \alpha (\varepsilon)}
$$
with $\, \alpha (\varepsilon) \, = \,  \log_{1/(r_i+ \varepsilon )} (\rho_i+ \varepsilon)$ and with $C(\varepsilon)>0$. Letting $\varepsilon \to 0$, we obtain the claim.

Next we establish the reverse inequality~$\alpha_{\varphi, J_i}
\le \log_{1/r_i} \, \rho_i$. Let $\varepsilon \in (0, r_i)$ and
$d_1,\ldots,d_k \in D(M)$, $k \in \n$. By Theorem~A2, there exist
$u \in U_i$ and a constant $C(u)>0$ such that  $\|A_{d_1} \cdots
A_{d_k}\, u\| \, \ge \, C(u) \rho_i^{\, k}$. Since the subspace
$U_i$ is spanned by the differences $v(y) - v(x), \ y - x \in
J_i$, $x,y \in G$, there exist $x_j, y_j \in G$, $y_j-x_j \in
J_i$, $j = 1, \ldots , n_i$ ($n_i$ dimension of $U_i$), such that
$\displaystyle u \, = \, \sum_{j=1}^{n_i} \gamma_j \bigl( v(y_j) -
v(x_j)\bigr)$, $\gamma_j \in \re$. Denote $x_j^{(k)} =
M_{d_1}^{-1}\cdots M_{d_k}^{-1} x_j\, $ and  $\, y_j^{(k)} =
M_{d_1}^{-1}\cdots M_{d_k}^{-1} y_j$. Thus, $\, x_j^{(k)},
y_j^{(k)} \in G_{d_1 \ldots d_k}$ and there exists
$C(\varepsilon)>0$ such that
\begin{equation} \label{aux300}
 \|y_j^{(k)} - x_j^{(k)}\| \, \le \, C(\varepsilon)\, (r_i-\varepsilon)^{-k}\|y_j - x_j\|, \quad k \in \n.
\end{equation}
Moreover, we have
\begin{eqnarray*}
&&\sum_{j=1}^{n_i} |\gamma_j|\, \cdot \, \bigl\|v(y_j^{(k)}) \, - \, v(x_j^{(k)})\bigr\| \  = \
\sum_{j=1}^{n_i} |\gamma_j| \, \cdot \, \bigl\|\, A_{d_1} \cdots A_{d_k}\, \bigl(v(y_j) \, - \, v(x_j)\bigr)\, \bigr\|  \ge \\
 &&\hspace{0.3cm}
\left\| \sum_{j=1}^{n_j} \gamma_j \, A_{d_1} \cdots A_{d_k} \, \bigl(v(y_j) \, - \, v(x_j)\,\bigr)\,  \right\| =
\left\|\, A_{d_1} \cdots A_{d_k}\, \left( \sum_{j=1}^n \gamma_j\, \bigl(v(y_j) \, - \, v(x_j)\,\bigr)\, \right)\,  \right\| =\\
&& \hspace{0.3cm}
\bigl\|\, A_{d_1} \cdots A_{d_k} u\, \bigr\|\ \ge \ C(u)\, \rho_i^{\, k}\,, \quad k \in \n\, .
\end{eqnarray*}
Consequently, at least one of the $n_i$ numbers $\|v(y_{j}^{(k)})
-  v(x_j^{(k)})\|, \ j = 1, \ldots , n_i$, is larger than or equal
to $\frac{C(u)}{\sum_j |\gamma_j|}\, \rho_i^{\, k}$. Combining
this estimate with \eqref{aux300} (i.e. $k \le -\frac{\log
\|y_j^{(k)} - x_j^{(k)}\|}{\log(r_i-\varepsilon)}+C_1$ with the
constant $C_1$ independent of $k$),
 we obtain
$$
\bigl\|v(y_{j}^{(k)})  \, - \,  v(x_j^{(k)})\| \ \ge \ C \, \bigl\|\, y_j^{(k)} \, - \, x_j^{(k)}\,
\bigr\|^{\, \alpha(\varepsilon)}\, ,
$$
where
$\alpha(\varepsilon) = \log_{1/(r_i-\varepsilon)}\, \rho_i$ and the  constant $C > 0$ does not depend on~$k$.
Since $\|y_j^{(k)} - x_j^{(k)}\|\to 0$ as $k \to \infty$, there exist arbitrary small segments~$[x_j^{(k)}\, , \, y_j^{(k)}]$
on which the variation of the function $v$ is at least a constant times the length of that segment to the
power of $\log_{1/(r_i-\varepsilon)}\, \rho_i$.
Therefore, $\alpha_{\varphi, J_i} \le \log_{1/(r_i-\varepsilon)} \, \rho_i$. Since $\varepsilon$ is arbitrary,
the claim follows.
\end{proof}

\begin{proof}[{\tt Proof of Theorem~\ref{th.holder}}]
We only show that the condition $\rho(\cA) < 1$ is sufficient for continuity of $\varphi$.
Let $\varepsilon \in (0,  1- \rho(\cA))$. Note that the set of rational $M$-adic points
$$
 Q=\bigcup_{k \in \n} Q_k, \quad Q_k=\left\{\sum_{j=1}^k M^{-j}d_j\ : \ d_j \in D(M)\right\}
$$
is dense in $G$. To determine the values of $v=v_\varphi$ in \eqref{eq.v} on $Q$ use the algorithm from subsection~\ref{ss.hold.u}.
We first show that~$v$ is uniformly bounded on~$Q$.
Denote  $C_0 = \max \{\|v(x) - v(y)\| \ :  \
x, y \in Q_1\}$. Then, for every $x = 0.d_1\ldots d_k \, \in Q_k \subset Q$, we have
\begin{eqnarray*}
&& \bigl\| \, v(0.d_1\ldots d_k \, - \, v(0)\, \bigr\| \ \le \
\sum_{j=1}^{k-1} \bigl\| T_{d_1}\cdots T_{d_j}
\bigl(\,  v(0.d_{j+1}) \, - \, v(0)  \, \bigr)   \, \bigr\|, \quad k \in \n.
\end{eqnarray*}
Note that, by construction, $v(0.d_{j+1})  -  v(0)\, = \, T_{d_{j+1}}v_0 - v_0 \, \in \, U$ for $j=1, \ldots, k-1$.
Therefore, by Theorem~A1, $\bigl\| T_{d_1}\cdots T_{d_j}
\bigl(v(0.d_{j+1})  -  v(0) \bigr)\bigr\|  \, \le \,
C_1 (\rho + \varepsilon)^j$, $j=1, \ldots, k-1$. Thus, we obtain
\begin{eqnarray*}
\bigl\| \, v(0.d_1\ldots d_k) \, - \, v(0)\, \bigr\| \ \le \
 C_1 \sum_{j=1}^{k-1}  (\rho(\cA) + \varepsilon)^j\ \le \
 C_1 \sum_{j=0}^{\infty} (\rho(\cA) + \varepsilon)^j \ = \ \frac{C_1}{1 - \rho(\cA) - \varepsilon}\, ,
\end{eqnarray*}
where  the constant $C_1>0$ is independent of $k$.
Hence, $\|v(x)\| \, \le \,
\|v(0)\|\, + \,  \frac{C_1}{1 - \rho(\cA) - \varepsilon}$
which proves  the uniform boundedness of $v$ on $G$.

The values of $v$ on $Q$ define the function $\varphi$ on
$\tilde Q$, where $\tilde{Q} = \cup_{k \in \z^s} (k+Q)$ is the set of all rational $M$-adic points of~$\re^s$.
The so constructed $\varphi: \tilde{Q} \rightarrow \re$ is supported on $K \cap \tilde Q$.
Using $\varphi$, define the extension $\tilde v: \tilde Q \to \re$ of $v$ in \eqref{eq.v}.
We show next that $\tilde v$ is uniformly continuous on $\tilde Q$, which implies that its extension to~$\re^s$
is continuous. Take arbitrary points $x, y \in \tilde Q$. By  Proposition~\ref{p.aux} and the same argument as in the
first part of the proof of Theorem~\ref{th.holder-direct}, we obtain
 $$
\bigl\|\tilde v(y) \, - \, \tilde v(x)\bigr\| \ \le \
 2 \, C_2 \, \|M\| \, \|v\|_{C(Q)}\, (\rho(\cA)+\varepsilon)^k\, ,
$$
where $k$ is the smallest number such that $\|M^k(y-x)\| \ge 1$.
Note that the value $\|v\|_{C(Q)}$ is finite because $v$ is bounded $Q$.
Since $\|M\|^k \ge 1/ \|y-x\|$, i.e. $k$ goes to $\infty$ as $\|y - x\|$ goes to zero,  $\tilde v$ is uniformly continuous on~$\tilde Q$,
which completes the proof of continuity. Thus, if $\rho(\cA) < 1$, then $\varphi \in C(\re^s)$.

By Theorem~\ref{th.holder-direct},
the Holder exponent $\alpha_\varphi$ of $\varphi$ on shifts along the subspace $J_i$ is equal to
 $\alpha_i = \log_{1/r_i}\, \rho_i$. We pass to a basis in the space $\re^s$,
 in which all the subspaces $J_i$ are orthogonal to each other.
Using  a natural expansion  $h = h_1 + \ldots + h_{q(M)}, \ h_i \in J_i$ we obtain for arbitrary
$\varepsilon > 0$
$$
\bigl\|\varphi (\cdot + h) - \varphi(\cdot)\bigr\| \ \le \ \sum_{i=1}^{q(M)} \, \bigl\|\varphi (\cdot + h_i) \, - \, \varphi(\cdot)\, \bigr\|\ \le \
\sum_{i=1}^{q(M)}  \, C \, \|h_i\|^{\, \alpha_i - \varepsilon} \ \le \ C \, q(M)^{\, \alpha - \varepsilon} \|h\|^{\, \alpha - \varepsilon},
$$
where $\alpha = \min\limits_{i=1, \ldots , q(M)}\alpha_i$. Consequently, $\alpha_{\varphi}\, = \, \min\limits_{i=1, \ldots , q(M)}\alpha_i$.
\end{proof}

\subsection{Examples} \label{ss.examples.C}

There are several types of anisotropic dilation matrices, e.g. hyperbolic dilations
or products of parabolic scaling and shear matrices, used in hyperbolic wavelet transform or 
discrete shearlet transform, see  \cite{Hyp_W, CoSauer, SauerK}. 

\begin{ex} {\rm In \cite{CoSauer}, the authors present a
general method for constructing interpolatory subdivision schemes with dilations which are  products of parabolic scaling and 
shear matrices.  We consider one of these examples with the anisotropic dilation matrix 
$M =\left(\begin{array}{cc} 2 &1 \\ 0 &3
\end{array}\right)$ and the refinement coefficients
\begin{eqnarray*}
 c(0,0) &=& \frac{1}{6}, \quad  c(0,1) = \frac{2}{6}, \quad  c(0,2) = \frac{3}{6}, \quad c(0,3) = \frac{4}{6},  \quad c(0,4) = \frac{5}{6}, \quad c(0,5) = 1, \\
 c(0,6) &=& \frac{5}{6}, \quad  c(0,7) = \frac{4}{6}, \quad  c(0,8) = \frac{3}{6}, \quad c(0,9) = \frac{2}{6},  \quad c(0,10) = \frac{1}{6}.
\end{eqnarray*} 
Note that the corresponding sequence of coefficients is supported on the set $\{(0, j) \ : \ j = 0, \ldots , 10\}$. By
\cite{GroH}, for the digit set $D(M)=\{(0,0), (1,0), (1, 1), (1, 2), (2, 1), (2, 2)\}$, the
set $G$ in \eqref{eq.G} is a tile. Using the result of
\cite{CHM2}, we determine the set $\Omega$ of size $19$. For the
corresponding $19\times 19$ transition matrices $T_d$, $d \in D(M)$, using Algorithm~1 from subsection \ref{ss.hold.u}
we construct the set $U$ starting with
$$
 v_0=\frac{1}{14} \left(\begin{array}{ccccccccccccccccccc}0&0&0&0&0&0&0&0&0&6&1&0&1&6&0&0&0&0&0 \end{array}\right)^T \in V,
 \quad ({\bf 1},v_0)=1.
$$
We obtain
\begin{eqnarray*}
 U^{(1)}&=& {\rm span}\{u_1^{(1)}, \ldots, u_4^{(1)} \}= {\rm span}\{T_{(1,1)}v_0-v_0, T_{(1,2)}v_0-v_0, T_{(2,1)}v_0-v_0, T_{(2,2)}v_0-v_0 
 \}, \\
 U^{(2)}&=&{\rm span} \{u_1^{(2)}, \ldots, u_{11}^{(1)}
 \}={\rm span} \{ U^{(1)}, T_{(0,0)}u_1^{(1)}, T_{(0,0)}u_2^{(1)}, T_{(0,0)}u_3^{(1)}, T_{(1,2)}u_1^{(1)}, T_{(2,1)} U^{(1)}\}, \\
  U^{(3)}&=& {\rm span} \{u_1^{(3)}, \ldots, u_{15}^{(3)}
 \}={\rm span}  \{ U^{(2)}, T_{(0,0)}u_6^{(2)},T_{(2,1)}u_8^{(2)},T_{(2,2)}u_9^{(2)}, T_{(2,2)}u_{10}^{(2)}\}, \\
  U=U^{(4)}&=&{\rm span}  \{u_1^{(4)}, \ldots, u_{16}^{(4)}
 \}={\rm span}  \{U^{(3)},
T_{(0,0)}u_{12}^{(3)}\}.
\end{eqnarray*}
Due to $\hbox{dim}(W) = 18$ and $\hbox{dim}(U)=16$,  $U$ is a proper subspace of $W$. Denote $\cA=\{A_d = T_d|U  \ : \ d \in D(M)\}$. We computed
the joint spectral radius of the set $\cA$ using the invariant polytope algorithm from~\cite{GP1}  and
obtained that the joint spectral radius is attained at the finite product of length $1$
$$
  \rho(\cA) \ = \ \rho(A_{(1,1)}) =0.75261 \dots 
$$
Since $M$ is anisotropic, there exist two non-zero subspaces $U_1, U_2 \subseteq U$. We construct these subspaces (each of dimension $15$) using Step $3$ of
Algorithm~1 from subsection \ref{ss.hold.u} with
$$ 
 U^{(1)} = {\rm span} \{v(0,0) - v(1/2, 0) \} \quad \hbox{or} \quad  U^{(1)} ={\rm span} \{v(0, 0)- v(1, 1)\},
$$
respectively. Note that, due to the continuity of $\varphi$, we have $v(0, 0) = v_0$. Furthermore,
$$
v(1/2, 0) = T_{(1,0)}v(0, 0), \quad  \hbox{due to}  \quad
M^{-1}  \left(\begin{array}{c} 1 \\ 0 \end{array}\right) = \left(\begin{array}{c} 1/2 \\ 0 \end{array}\right) \in G,
$$
and
$$
v(1,1) =  T_{(2,2)}v(1,1), \quad  \hbox{due to} \quad
 \sum_{j=1}^\infty M^{-j}  \left(\begin{array}{c} 2 \\ 2 \end{array}\right) = \left(\begin{array}{c} 1 \\ 1 \end{array}\right) \in G.
$$
For the restrictions $\cA|_{U_1}$ and $\cA|_{U_2}$, we obtain $\rho_1 = \rho_2 =  \rho(\cA)$. Therefore, 
by Theorem~\ref{th.holder},
$$
\alpha_{\varphi} \ =  \log_{1/3} \rho(A_{1,1}) \ = \ 0.25869 \ldots
$$}
\end{ex}

\noindent The next example shows that in some cases $W=U=U_1=U_2$.

\begin{ex} {\rm
 We consider the dilation matrix $M=\left(\begin{array}{rr}2&1\\
1&-1\end{array}\right)$ and the refinement coefficients
$$
c_{(0,0)}=\frac{1}{2}, \quad  c_{(1,-1)}=\frac{1}{2}, \quad  c_{(1,0)}=1, \quad
 c_{(2,0)}=\frac{1}{4} \quad \hbox{and} \quad c_{(1,1)}=\frac{3}{4}.
$$
The dilation matrix has eigenvalues $\lambda_1=\frac{1 -
\sqrt{13}}{2}$ and $\lambda_2=\frac{1 + \sqrt{13}}{2}$. By
\cite{GroH}, for the digit set $D(M)=\{(0,0), (1,0), (2,0)\}$, the
set $G$ in \eqref{eq.G} is a tile. Using the result of
\cite{CHM2}, we determine
$$
 \Omega=\{-1,0,1\}^2 \cup \{(0,\pm 2),(1,-2),(-1,2)\}.
$$
The corresponding transition matrices $T_{(0,0)}$, $T_{(1,0)}$ and
$T_{(2,0)}$ are given by\scriptsize
$$
 \frac{1}{4}\, \left(\begin{array}{ccccccccccccc}
  2&0&0&0&0&0&0&0&0&2&4&3&0\\
  0&3&1&0&0&0&4&0&0&0&0&0&0\\
  2&0&0&0&4&2&0&3&0&0&0&1&0\\
  0&0&0&0&0&2&0&1&4&0&0&0&0\\
  0&0&0&0&0&0&0&0&0&2&0&0&0\\
  0&0&3&4&0&0&0&0&0&0&0&0&1\\
  0&1&0&0&0&0&0&0&0&0&0&0&0\\
  0&0&0&0&0&0&0&0&0&0&0&0&0\\
  0&0&0&0&0&0&0&0&0&0&0&0&3\\
  0&0&0&0&0&0&0&0&0&0&0&0&0\\
  0&0&0&0&0&0&0&0&0&0&0&0&0\\
  0&0&0&0&0&0&0&0&0&0&0&0&0\\
  0&0&0&0&0&0&0&0&0&0&0&0&0 \end{array} \right), \quad
 \frac{1}{4}\, \left(\begin{array}{ccccccccccccc}
  4&2&2&0&3 &0&2&0&0&0 &1&0&0\\
  0&0&0&0&0 &0&0&0&0&0 &0&0&0\\
  0&2&0&0&1 &4&2&0&3&0 &0&0&0\\
  0&0&0&0&0 &0&0&0&1&0 &0&0&0\\
  0&0&2&0&0 &0&0&0&0&4 &3&0&2 \\
  0&0&0&1&0 &0&0&0&0&0 &0&0&0\\
  0&0&0&0&0 &0&0&0&0&0 &0&0&0\\
  0&0&0&0&0 &0&0&0&0&0 &0&0&2 \\
  0&0&0&3&0 &0&0&0&0&0 &0&0&0\\
  0&0&0&0&0 &0&0&2&0&0 &0&2&0 \\
  0&0&0&0&0 &0&0&0&0&0 &0&2&0\\
  0&0&0&0&0 &0&0&0&0&0 &0&0&0\\
  0&0&0&0&0 &0&0&2&0&0 &0&0&0 \end{array} \right)
$$
\normalsize and \scriptsize
$$
 \frac{1}{4}\, \left(\begin{array}{ccccccccccccc}
  1&4&0&0&0 &0&4&0&0&0 &0&0&0\\
  0&0&0&0&0 &0&0&0&0&0 &0&0&0\\
  0&0&0&0&0 &0&0&0&0&0 &0&0&0\\
  0&0&0&0&0 &0&0&0&0&0&0&0&0\\
  3&0&4&2&0 &0&0&0&0&1 &0&0&0 \\
  0&0&0&0&0 &0&0&0&0&0 &0&0&0\\
  0&0&0&0&0 &0&0&0&0&0&0&0&0\\
  0&0&0&2&0 &0&0&0&0&3 &0&0&4 \\
  0&0&0&0&0 &0&0&0&0&0 &0&0&0\\
  0&0&0&0&2 &0&0&4&2&0 &0&0&0 \\
  0&0&0&0&2 &4&0&0&0&0 &2&4&0\\
  0&0&0&0&0 &0&0&0&0&0 &2&0&0\\
  0&0&0&0&0 &0&0&0&2&0 &0&0&0 \end{array} \right),
$$
\normalsize respectively. The matrix $T_{(0,0)}$ has one
eigenvalue $1$ with the corresponding eigenvector
$$
 v_0=\frac{1}{5} \left(\begin{array}{ccccccccccccc}0&4&0&0&0&0&1&0&0&0&0&0&0 \end{array}\right)^T \in V,
 \quad ({\bf 1},v_0)=1.
$$
Using Algorithm~1 from subsection
\ref{ss.hold.u}, we obtain $U$ with $\hbox{dim}(U)=12$. Thus, $W=U$ and, therefore, 
$\cA=\{A_d=T_d|_{W} \ : \ d \in D(M)\}$. 
We computed the joint spectral radius of the set $\cA$ using 
the invariant polytope algorithm from~\cite{GP1} and obtained that 
the joint spectral radius is attained at the finite product 
$(A_{(0,0)} A_{(1,0)})^2 A_{(0,0)}^2 A_{(2,0)}$ of length~$7$, i.e.
$$
   \rho(\cA) \ = \  \Bigl[ \rho \, \bigl( (A_{(0,0)} A_{(1,0)})^2 A_{(0,0)}^2 A_{(2,0)}\bigr) \Bigr]^{1/7} \ = \ 
   0.93816\ldots .
$$
We verified  that the matrix family $\cA$ is irreducible in this case. Hence, the only non-zero common invariant 
subspace of the matrices in $\cA$ is $U$. Thus, $U_1 = U_2 = U = W$. Therefore, $\rho_1 = \rho_2 =  \rho(\cA)$ and
Theorem 1 implies that 
$$
\alpha_{\varphi} \ = \ \log_{1/\rho(M)}\, \rho(\cA) \ = \ 
\frac17\, \log_{\frac{2}{1 + \sqrt{13}}} \rho \, \bigl( (A_{(0,0)} A_{(1,0)})^2 A_{(0,0)}^2 A_{(2,0)}\bigr) \ = \ 
0.07652\ldots 
$$ }
\end{ex}

\section{Higher order regularity}\label{s.higher}

\noindent In this section, we show that the derivatives of the multivariate refinable function~$\varphi \in \cS'(\re^s)$
satisfy a system of nonhomogeneous refinement equations. The differentiability of~$\varphi \in C(\re^s)$ is then equivalent to  continuity of the solutions of all these equations, see Theorem \ref{th.deriv}. The main idea is
that the directional derivatives of $\varphi$
along the eigenvectors of the dilation matrix $M$ satisfy certain refinement equations and the directional derivatives
along the generalized eigenvectors (of the Jordan basis) of $M$ satisfy nonhomogeneous refinement equations, see
Proposition~\ref{p.derivatives}.

\begin{defi}\label{d.gen-ref}
A multivariate nonhomogeneous refinement equation is a functional equation of the form
$$
\varphi \, = \, \bT \, \varphi\, + \, g,
$$
where $\bT$ is the transition operator in~(\ref{eq.transit}) and $g$ is a compactly supported function or distribution.
\end{defi}

\noindent For more details on nonhomogeneous refinement equations see e.g. \cite{DH,JJS, JiangLi, SZ} and references therein.

\noindent  Let $
E=\{e_1, \ldots, e_s\} $
be the Jordan basis of the matrix~$M$ in~$\re^s$. The Jordan basis consists of
the eigenvectors of $M$, which satisfy $Me_i = \lambda e_i$, and of the generalized eigenvectors, which
satisfy $Me_i = \lambda e_i + e_{i-1}$. Consider an $\ell \times \ell$ Jordan block of $M$ corresponding to an eigenvalue $\lambda$. With a slight abuse of notation, denote by
$$
 E_\lambda=\{e_{1}, \ldots, e_{\ell}\} \subset E, \quad Me_{1} = \lambda e_{1}, \quad
 Me_{i} = \lambda e_{i} + e_{i-1}, \ i = 2, \ldots , \ell,
$$
the Jordan basis corresponding to this Jordan block. In the following we study the properties of the
directional derivatives of the refinable function $\varphi \in \cS'(\re^s)$, which belong to the following
subspaces of $\cS'(\re^s)$.

\begin{defi} For a vector $a \in \re^s$, we denote by
$$
\cS'_a(\re^s)\ =\ \left\{ \varphi \in \cS'(\re^s) \ : \ \int_{x=at+b, \ t\in \re^s} \varphi(x)dx=0, \ b \in \re^s \right\}
$$
the space of compactly supported distributions, whose mean along every straight line
$\bigl\{x=at+b \ : \ t \in \re\bigr\}$, $\ b\in \re^s$,  parallel to~$a$, is equal to zero.
\end{defi}

\noindent By $\nabla \varphi\, = \, \bigl(\frac{\partial \varphi}{\partial x_1}, \ldots ,
\frac{\partial \varphi}{\partial x_s} \bigr)$ we denote the total derivative (gradient) of $\varphi$ and by
$\frac{\partial \, \varphi}{\partial \, a}\, = \, \bigl( a\, , \, \nabla \varphi \bigr)$, its directional
derivative along a nonzero vector $a \in \re^s$. Due to the compact support of $\varphi \in \cS'(\re^s)$, its directional derivative
$\frac{\partial \, \varphi}{\partial \, a}$ belongs to  $\cS'_a(\re^s)$.

\noindent The next result shows that a directional derivative of a refinable function~$\varphi$ along an eigenvector of the
dilation matrix $M$ is also a refinable function and satisfies the refinement equation~(\ref{eq.e1}).
A directional derivative of~$\varphi$ along a generalized eigenvector of $M$ satisfies the nonhomogeneous
refinement equation~(\ref{eq.ei}).

\begin{prop}\label{p.derivatives} Let $\varphi \in \cS'(\re^s)$ and $\lambda$ be an eigenvalue of $M$.
If $\varphi=\bT \varphi$, then $\varphi_i = \bigl(e_{i}, \nabla \varphi\bigr) \in \cS'_i:=\cS_{e_{i}}'(\re^s)$, $e_{i} \in E_\lambda$,
satisfy the refinement equation
\begin{equation}\label{eq.e1}
\varphi_1 \ = \ \lambda \, \bT \varphi_1\,
\end{equation}
and the nonhomogeneous refinement equations
\begin{equation}\label{eq.ei}
\varphi_i \ = \ \lambda \, \bT \varphi_i \ + \  \sum_{j=1}^{i-1}(-1)^{j-1}\, \lambda^{-j}\, \varphi_{i-j}\, ,
\quad i = 2, \ldots , \ell\ =\ { \rm dim}\, (E_\lambda).
\end{equation}
Conversely, the system of equations~(\ref{eq.e1})-(\ref{eq.ei}) possesses a unique up to a normalization solution
$(\varphi_1, \ldots, \varphi_\ell) \in \cS'_1 \times \ldots \times \cS'_\ell$. Moreover, if $\varphi \in \cS'(\re^s)$ satisfies $\bigl(e_{i}, \nabla \varphi\bigr)=\varphi_i$, $e_i \in E_\lambda$, then $\varphi=\bT \varphi$ along the lines parallel to
$e_i \in E_\lambda$.
\end{prop}
\begin{proof} By induction on $\ell$, we show that, if $\varphi \in \cS'(\re^s)$ satisfies the refinement equation
$\varphi=\bT \varphi$, then $\varphi_i = \bigl(e_{i}, \nabla \varphi\bigr)$, $i=1,\ldots,\ell$, satisfy (\ref{eq.e1})-(\ref{eq.ei}).
For $\ell=1$, due to $Me_{1} = \lambda e_{1}$, we have
\begin{eqnarray*}
\varphi_1&=& \bigl(e_{1}, \nabla \varphi\bigr) \ = \
\left( e_{1}, \sum_{k \in\z^s}  c_k M^T \nabla \varphi \right)(M\cdot - k) =
\left( M\, e_{1}, \sum_{k \in \z^s} c_k \nabla \varphi(M\cdot - k)\, \right) \\ &=& \
\lambda \, \sum_{k\in \z^s} \, c_k\, \Bigl( \, e_{1}, \nabla \varphi\, \Bigr)(M\cdot - k) \ = \
\lambda\, \sum_{k \in\z^s} c_k \varphi_1 (M \cdot - k)\, =\lambda \bT \varphi_1\, .
\end{eqnarray*}
For $i \in \{2, \ldots,\ell\}$, due to  $Me_{i} = \lambda e_{i} + e_{i-1}$ and \eqref{eq.ei},
we obtain
\begin{eqnarray*}
 \varphi_i \, = \, \lambda \bT \varphi_{i} \, + \, \bT \varphi_{i-1}=\lambda \bT \varphi_{i} \, + \lambda^{-1}\varphi_{i-1}-\sum_{j=1}^{i-2} (-1)^{j-1}\lambda^{-j} \varphi_{i-1-j}.
\end{eqnarray*}

Conversely, by \cite{JJS}, the system~(\ref{eq.e1})-(\ref{eq.ei}) possess a unique up to normalization solution. We show that the compactly supported primitive~$\varphi$ of $\varphi_i$, $i=1,\ldots,\ell$, along $e_{i} \in E_\lambda$ satisfies~$\varphi = \bT \varphi$. Indeed, since $Me_{1} = \lambda e_{1}$ and by (\ref{eq.e1}), we obtain
$$
\left( e_1 ,  \nabla \varphi \right) \ = \lambda \bT
\left( e_1, \nabla \varphi \right) \ = \left( e_1\, , M^T \bT \nabla \varphi \right)=
\left( e_1\, , \nabla \bT \varphi \right),
$$
which implies that the gradient of the function $\bT \varphi - \varphi$ is orthogonal to $e_{1}$.
Hence, the function $\bT \varphi - \varphi$  is constant along all lines parallel to $e_{1} \in E_\lambda$.
On the other hand, $\varphi$ is compactly supported, consequently $\bT \varphi - \varphi=0$ along all lines parallel
to $e_{1} \in E_\lambda$. For $i \in \{2, \ldots,\ell\}$, due to $Me_{i} = \lambda e_{i}+ e_{i-1}$ and
$\varphi_i=\lambda \bT \varphi_i+\bT \varphi_{i-1}$, we obtain
$$
 \left(e_{i},  \nabla \varphi \right) \ = \lambda \bT
\left(e_{i}, \nabla \varphi \right)+  \bT\left(e_{i-1}, \nabla \varphi \right)
= \left( \lambda e_i+e_{i-1}\,  \bT\nabla \varphi \right)=
\left( e_i\, , \nabla \bT \varphi \right).
$$
Analogous considerations about $\bT \varphi-\varphi$ along $e_i \in E_\lambda$ imply the claim.
\end{proof}

\begin{remark}\label{r.250} {\em
If, for the eigenvalue $\lambda$ of the dilation matrix $M$, the set $E_\lambda$ does not
contain any generalized eigenvectors, then
the system (\ref{eq.e1})-(\ref{eq.ei}) reduces to homogeneous refinement equations
$\varphi_i \ = \ \lambda \, \bT \varphi_i$,  $e_i \in E_\lambda$. }
\end{remark}

\noindent The main result of this section, Theorem \ref{th.deriv}, states that $\varphi \in C^1(\re^s)$ if and only if
the (nonhomogeneous) refinement equations in Proposition \ref{p.derivatives} corresponding to the Jordan basis $E \subset \re^s$ of the
dilation matrix $M$ have continuous solutions $\varphi_i \in \cS'_{e_i}(\re^s)$, $e_i \in E$.
The directional derivatives $\varphi_i$, $i=1,\ldots,s$, determine the total derivative $\nabla \varphi$
of $\varphi$. Moreover, $\varphi_i$, $i=1,\ldots,s$, can be constructed and their H\"older exponents can be computed
as described in section~\ref{s.hold} (see Remark~\ref{r.610}). Thus, the higher regularity of any refinable function $\varphi$ can be analyzed
by this recursive reduction to a set of continuous refinable functions.

\begin{theorem}\label{th.deriv} Let $\varphi \in \cS'(\re^s)$.
There exist continuous solutions $\varphi_i \in \cS'_{e_i} (\re^s), \, e_i \in E_\lambda$, of
(\ref{eq.e1}) -- (\ref{eq.ei}) for each eigenvalue $\lambda$ of the dilation matrix $M$ if and only if $\varphi \in C^1(\re^s)$
satisfies $\varphi=\bT \varphi$ and
 $\frac{\partial \, \varphi}{\partial \, e_i}\, = \, \varphi_i$, $e_i \in E$.
\end{theorem}
\begin{proof}
If $\varphi \in C^1(\re^s)$  is a compactly supported solution of the refinement equation
$\varphi=\bT \varphi$, then, by Proposition~\ref{p.derivatives}, its directional derivatives $\varphi_i:=\frac{\partial \, \varphi}{\partial \, e_i}
\in C(\re^s)$ along $e_i \in E$ satisfy equations~(\ref{eq.e1})-(\ref{eq.ei}).
 Conversely, if the equations~(\ref{eq.e1})-(\ref{eq.ei}) possess continuous solutions,
then, by Proposition~\ref{p.derivatives}, $\varphi$ is in $C^1(\re^s)$, $\frac{\partial \, \varphi}{\partial \, e_i}=\varphi_i$, $e_i \in E$, and satisfies $\varphi=\bT \varphi$.
\end{proof}

\begin{cor}\label{c.deriv}
Suppose that $E$ does not contain any generalized eigenvectors, i.e., the matrix $M$
has a basis of eigenvectors; then
\begin{description}
\item[$(i)$] If $\varphi \in C^1(\re^s)$ is refinable, then $\varphi_i \, = \, \frac{\partial\, \varphi}{\partial \, e_i}\in C(\re^s)$ satisfy $\varphi_i \ = \ \lambda_i \, \bT \varphi_i$, $i=1, \ldots,s$.
\item[$(ii)$] Conversely, if the solutions $\varphi_i \in \cS_{e_i}'$ of $\varphi_i \ = \ \lambda_i \, \bT \varphi_i$,
$i=1, \ldots,s$, are continuous, then the solution of~$\varphi =\bT \varphi$
belongs to $C^1(\re^s)$. Moreover, $\frac{\partial \varphi}{\partial e_i}=\varphi_i$, $i=1, \ldots,s$.
\end{description}
\end{cor}

\begin{remark}\label{r.610}
{\em The system of refinement equations~(\ref{eq.e1})-(\ref{eq.ei}) is solved and analysed in the same
way described in subsection~\ref{ss.hold.u}.
First we solve the equation $\varphi_1 = \lambda \bT \varphi_1$. We find $v_{\varphi_1}(0)$
as an eigenvector of the matrix $T_0$ with the eigenvalue~$1/\lambda $. If
$T_0$ does not have this eigenvalue, then  equation $\varphi_1 = \lambda \bT \varphi_1$ does not have a solution, and hence $\varphi \notin C^1(\re^s)$.
Then we define the space $U_{\lambda , 1}$ as the minimal common invariant subspace of $\cT$ containing the vectors $\lambda T_d  v_{\varphi_1}(0)\, - \, v_{\varphi_1}(0), \, d \in D(M)$. This can be done by Algorithm~1
from subsection~\ref{ss.hold.u}. Then $\varphi_1 \in C(\re^s)$ if and only if
the joint spectral radius of the matrices $\lambda T_d,\, d \in D(M),$
restricted to the subspace~$U_{\lambda , 1}$ is smaller than one. The H\"older regularity
of $\varphi_1$ is computed by formula~(\ref{eq.holder}) for the matrices $A_d = \lambda T_d|_{U_{\lambda , 1}}$.
Similarly, we solve the other equations of the system~(\ref{eq.ei}) successively  for $i = 2, \ldots , \ell$. }
\end{remark}

\section{Existence and smoothness  in~$\mathbf{L_p(\re^s),\ 1 \le p < \infty}$}\label{s.p}

\noindent In this section, see Theorem \ref{th.p}, we characterize the existence of refinable functions in
$\varphi \in L_p(\re^s)$, $1 \le p < \infty$, and provide a formula for the
H\"older exponent of such $\varphi$, see Theorem~\ref{th.holder-p}, in terms of the $p$-radius
($p$-norm joint spectral radius  \cite{J95, P97}) of a set of transition matrices.

\begin{defi}\label{d.pr}
    For $1 \le p < \infty$, the $p$-radius ($p$-norm joint spectral radius)
    of a finite family of linear
    operators~$\cA\, = \, \{A_0, \ldots , A_{m-1}\}$ is defined by
    $$
    \rho_p=\rho_p(\cA)\ = \ \lim_{k \to \infty}\, \Bigl(\,
    m^{-k}\sum_{A_{d_i} \in \cA, \, i = 1, \ldots , k}\|A_{d_1}\cdots A_{d_k}\|^p \, \Bigr)^{\, 1/pk}.
    $$
\end{defi}

\noindent Note that, for $\varphi\in L_p(\re^s)$, the difference space $U$ is defined similarly to \eqref{eq.U} by

\begin{equation}\label{eq.U_Lp}
 U\ = \ {\rm span}\, \Bigl\{\, v(y)\, - \, v(x) \quad :  \quad \hbox{for almost all} \quad y, x \, \in G\, \Bigr\} \subseteq W.
\end{equation}

\noindent Although the estimates in Theorem~\ref{th.holder-p} look similar to the ones from section \ref{s.hold}, the corresponding proofs require totally different techniques.

\noindent The H\"older exponent of a function $\varphi \in L_p(\re^s)$  is defined by
$$
 \alpha_{\varphi, p}\quad = \quad \sup \Bigl\{\, \alpha \ge 0 \quad : \quad \bigl\|\varphi(\cdot  + h) \, - \, \varphi(\cdot) \bigr\|_p \ \le \  C\,  \|h\|^{\, \alpha}\, , \quad h \in \re^s\, \Bigr\}\, .
$$
We use the notation~$\|\cdot \|_{L_p(\re^s)}\, = \,
\|\cdot \|_{p}$.
To determine the influence of the dilation matrix $M$ on the  H\"older exponent of  $\varphi \in L_p(\re^s)$,  in Theorem \ref{th.holder-p},
we consider the H\"older exponents of $\varphi$ along the subspaces determined by the Jordan basis of $M$.
The H\"older exponent of $\varphi$ along a subspace $J \subset \re^s$ is defined by
$$
\alpha_{\varphi, J, p}\quad = \quad \sup
\Bigl\{\, \alpha \ge 0 \quad : \quad \bigl\|\varphi(\cdot  + h) \, - \, \varphi(\cdot)\bigr\|_p \ \le \  C\,  \|h\|^{\, \alpha}\, , \quad h \in J\, \Bigr\}\, .
$$

\noindent In the proofs of Theorems \ref{th.p} and \ref{th.holder-p} we use the following auxiliary results.

\subsection{Auxiliary results}\label{ss.p.prem}

\noindent The following analogues of Theorems~A1 and A2 from section \ref{s.hold} were proved in~\cite{P08}.

\noindent \textbf{Theorem A3}. {\it Let $1 \le p < \infty$. For a finite family~$\cA$ of $m$ operators acting in~$\re^n$
and for any $\varepsilon > 0$, there exists a norm $\|\cdot \|_{\varepsilon}$ in $\re^n$ such that
$$
\Bigl( m^{-1}\,  \sum_{A \in \cA}\, \|A u\|_{\varepsilon}^p \, \Bigr)^{\, 1/p}\ \le \
 \bigl(\rho_p + \varepsilon\bigr)\, \|u\|_{\varepsilon}\, , \qquad u \in \re^n\, .
$$
}
\smallskip

\noindent For~$1 \le p < \infty$ and for a finite family~$\cA$ of $m$  operators acting in~$\re^n$, we denote
$$
\cF_{p} (k,u)\quad = \quad \Bigl(\,
m^{-k}\, \sum_{ A_{d_i} \in \cA, i=1,\ldots,k}
\, \bigl\|A_{d_1}\cdots A_{d_k}u\bigr\|^p \ \Bigr)^{1/p} , \qquad k \in \n\, .
$$
Since each norm~$\|\cdot\|$ in $\re^n$ is equivalent to the norm $\|\cdot \|_{\varepsilon}$,
Theorem~A3 yields the following result.

\begin{cor}\label{c.p-rad} Let~$1 \le p < \infty$.
For every $\varepsilon > 0$, there exists a constant $C(\varepsilon) >0$ such that
$\cF_{p} (k,u) \, \le \, C(\varepsilon) (\rho_p + \varepsilon)^k\|u\|$ for all
$u \in \re^n$ and $k \in \n$.
\end{cor}

\noindent \textbf{Theorem A4}. {\it Let~$1 \le p < \infty$ and $\cA$ be a finite family of $m$ linear operators in~$\re^n$.
Then for every $u\in \re^n$ that does not belong to any common
invariant  linear subspace of $\cA$, there exists a constant $C(u) > 0$ such that
\begin{equation}\label{eq.constant}
 \cF_p(k, u) \ \ge \ C(u)\, \rho_p^{\, k}, \qquad  k \in \n.
\end{equation}
 }

\noindent In Lemma \ref{l.p-rad} we relax assumptions of Theorem A4.
We show that~(\ref{eq.constant}) holds for all points $u$ apart from the ones in a proper linear subspace of $\re^n$.

\begin{lemma}\label{l.p-rad}  Let~$1 \le p < \infty$.
Every finite family $\cA$ of $m$ linear operators possesses a common invariant linear subspace
$\cL \subset \re^n$, $\cL \not =\re^n$, (possibly $\cL = \{0\}$) such that for every $u \notin \cL$ there
exists a constant  $C(u) > 0$ for which~(\ref{eq.constant}) holds.
\end{lemma}
\begin{proof} Without loss of generality, after a suitable normalization, it
can be assumed that $\rho_p = 1$. Let $\cL$ be the biggest by inclusion common invariant subspace
of $\cA$ such that $\rho_p(\cA|_{\cL}) < 1$. Note that $\cL$ is a proper subspace of $\re^n$, otherwise
we get a contradiction to $\rho_p = 1$. Hence, ${\rm dim}\, \cL \, \le \, n-1$.
Take  arbitrary $u \notin \cL$ and denote by $\cL_u$ the minimal common invariant subspace of
$\cA$ that contains $u$. If $\rho_p(\cA|_{\cL_u}) < 1$, then the $p$-joint spectral radius of $\cA$
on the linear span of $\cL$ and $\cL_u$ is equal to $\max\, \{\rho_p(\cA|_{\cL}), \rho_p(\cA|_{\cL_u})\} < 1$,
which contradicts the maximality of~$\cL$. Hence, $\rho_p(\cA|_{\cL_u}) = 1$.
Since $u$ does not belong to any common invariant subspace of the finite family $\cA|_{\cL_u}$,
by Theorem A4, there exists a constant $C(u)>0$ such that
$\cF_p(k, u) \ge  C(u) \rho(\cA|_{\cL_u})^k, \, k \in \n$. On the other hand, $\rho_p=\rho(\cA|_{\cL_u}) = 1$,
and, hence, the claim follows.
\end{proof}

\noindent For~$1 \le p < \infty$, in the rest of this section, we denote  by $C(u) > 0$ the largest possible constant in inequality~(\ref{eq.constant}), i.e.,
\begin{equation} \label{d.C}
C(u)\, = \, \inf_{k \in \n}\, (\rho_p)^{\, -k}\cF_p(k, u).
\end{equation}
This function is upper semi-continuous and, therefore, is measurable.

\subsubsection{Properties of the space $U$} \label{ss.p.U}

Note that, by~(\ref{eq.ss1}), the vector $\displaystyle z = \int_{G}v(x)\, dx \in V$ is an eigenvector of the operator
$\displaystyle T \, = \, \frac{1}{m}\, \sum_{d \in D(M)}\, T_d  $ associated with the eigenvalue~$1$.
The following result  is an analogue of Proposition~\ref{p.U0}.

\begin{prop}\label{p.U0p}
If $\varphi \in L_p(\re^s)$, $1 \le p < \infty$, then the subspace $U$ in \eqref{eq.U_Lp} coincides with the smallest by
inclusion common invariant subspace of the matrices $T_d, \, d \in D(M)$, and contains the $m$ vectors $\, T_dz - z, \, d \in D(M)$,
where the nonzero $z \in V$ satisfies $Tz=z$.
\end{prop}

\noindent  The proof of Proposition \ref{p.U0p} is similar to the one of Proposition~\ref{p.U0}. Note that,
since $\displaystyle \sum_{k \in \z^s} c_k = m$, the column sums of the matrix $T$ are equal to one. Hence, $T$
has at least one eigenvalue one. This eigenvalue does not have to be simple. Nevertheless, Proposition \ref{p.U3p}
guarantees that $U$ is well-defined.  The proof of Proposition \ref{p.U3p}  is similar to the one of
Proposition~\ref{p.U3}.

\begin{prop}\label{p.U3p}
There exists at most one eigenvector $z \in V$ of $\, T$ associated to the eigenvalue~$1$ such that~$\rho_p(\cA) < 1$.
\end{prop}

\subsection{$\mathbf{L_p}$-solutions of refinement equations}\label{ss.p.Lp}

\begin{theorem}\label{th.p}
A refinable function $\varphi$ belongs to~$L_p(\re^s)$, $1 \le p < \infty$, if and only if $\rho_p(\cA) < 1$.
\end{theorem}
\begin{proof} Assume first that $\rho_p=\rho_p(\cA) < 1$. Choose $\varepsilon \in (0,  1- \rho_p)$  and consider the norm $\|\cdot\|_{\varepsilon}$ in $U$ as in Theorem~A4. Define the function space
$$
V_{U, p}\ = \
\bigl\{f \in L_p(G+\Omega) \ : \ v_f(x) \in V, \ v_f(x) - v_f(y) \in U\quad \mbox{ a.e.}, \ x, y \in G\, \bigr\}
$$
with the norm
$\|f\| \, = \, \bigl(\, \int_{G}\|v_f(x)\|^p_{\varepsilon}\, d x\, \bigr)^{1/p}$.
The space $V_{U, p}$ is nonempty because it at least contains a piecewise constant function $f$
such that $v_f \equiv z$ a.e., where $z \in V$ is the eigenvector of the
operator $\displaystyle \frac{1}{m}\sum_{d \in D(M)}\, T_d$ associated with the eigenvalue one.
Note that, for $f_1, f_2 \in V_{U, p}$, we have
$$
 \|\bT (f_1 - f_2)\| \, = \, \left( \int_{G} \|\bA (v_{f_1} - v_{f_2})(x)\|_{\varepsilon}dx \right)^{1/p} \, \le \,
(\rho_p + \varepsilon)\|f_1-f_2\|.
$$
Therefore, due to $\rho_p + \varepsilon<1$, $\bT$ is a contraction on $V_{U, p}$,  and, hence, it has a unique fixed point $\varphi$,
which is the solution of the refinement equation~$\varphi=\bT\varphi$.

Assume next that $\varphi \in L_p(\re^s)$.
By Lemma~\ref{l.p-rad}, there exists a proper subspace
$\cL \subset U$ (due to $\hbox{dim} U=n$, we associate $\re^n$ with $U$) invariant under $\cA$ such that
$\cF_p(k, u) \ge C(u)\, \rho_{p}^{\, k}, \ k \in \n$, whenever
$u \notin \cL$. Since $U$ is the smallest by inclusion subspace
of $\re^N$ invariant under $\cA$ and containing the differences $v(x_2) - v(x_1)$ for almost all
$x_1, x_2 \in G$, the set
$
 \{(x_1, x_2) \in G^2 \ : \  v(x_2) - v(x_1) \notin \cL \}$
has a positive Lebesque measure $\mu$ in~$\re^s \times \re^s$.
Hence, the set
\begin{equation} \label{aux71}
 \{(x, h) \in \re^s\times \re^s \ : \ x, x+h \in G, \ v(x+h) - v(x) \notin \cL\}
\end{equation}
has a positive Lebesque measure. By the Fubini theorem, the set in \eqref{aux71} has sections of positive Lebesque
measure. Thus, there exists $h \in \re^s$ such that
$$
\mu \bigl\{x \in G \ : \ x + h \in G \, , \  v(x+h) - v(x) \notin \cL\bigr\} >0.
$$
Therefore, there exist $\delta > 0$ and a set $H \subset G$ of positive Lebesque measure
such that the function in \eqref{d.C} satisfies $C\bigl( v(x+h) - v(x)\bigr) > \delta$ for almost all~$x \in H$.
Thus,
\begin{equation}\label{eq.p-low}
\cF_p\, \bigl(k, v(x+h) - v(x)\bigr) \quad \ge \quad
\delta \,  \rho_{p}^k, \qquad k \in \n \ ,\quad  \mbox{for almost all} \ x \in H\, .
\end{equation}
Denote $h_k = M^{-k}h$, $k \in \n$, then, by (\ref{eq.ss1}), we get
\begin{eqnarray}\label{eq.p1}
\bigl\| \,  v(\cdot+h_k) \, - \,   v\, \bigr\|^p_{L_p(\re^s)} \ &\ge&
\bigl\| \,  v(\cdot+h_k) \, - \,   v\, \bigr\|^p_{L_p(G)} \notag \\ &\ge& \
\sum_{d_1, \ldots, d_k \in D(M)}\int_{H_{d_1 \ldots d_k}}\, \bigl\| \,  v(x+h_k) \, - \,   v(x)\, \bigr\|^{p}\, dx\,
\notag \\
&=&\sum_{d_1, \ldots, d_k \in D(M)}m^{-k}\, \int_{H}\, \bigl\|A_{d_1}\cdots A_{d_k}\,
\bigl( \,  v(y+h) \, - \,   v(y)\, \bigr) \bigr\|^p\, dy \notag\\& =& \
\int_{H}\, \cF^p_p\, \bigl( k\, , \, \bigl( \,  v(y+h) \, - \,   v(y)\, \bigr)\, \bigr)\, dy\, .
\end{eqnarray}
Since $y \in H$, by \eqref{eq.p-low}, we obtain
$$
\int_{H}\, \cF^p_p\, \bigl( k\, , \, \bigl( \,  v(y+h) \, - \,   v(y)\, \bigr)\, \bigr)\, dy\ \ge \
\delta^p\, \rho_p^{ p\, k}\, \mu(H)\, .
$$
Thus,
\begin{equation}\label{eq.p2}
\bigl\| \,  v(\cdot+h_k) \, - \,   v\, \bigr\|_{L_p(\re^s)} \quad \ge \quad
\delta\ \rho_p^{k}\ [\mu(H)]^{1/p}\ , \qquad h_k = M^{-k}h\, , \quad k \in \n\, .
\end{equation}
Since $v \in L_p(\re^s, \re^N)$ and $h_k$ goes to $0$ as $k \to \infty$, we get $\rho_p < 1$.
\end{proof}

\begin{remark}\label{r.30}
{\em The proof of Theorem \ref{th.p} is much simpler than that of Theorem~\ref{th.holder}.
Indeed, the argument with a contraction operator $\bT$ on the affine subspace $V(U, p)$
cannot be directly extended to prove the continuity of $\varphi$ due to the following reason:
the piecewise constant function $f$ for which $v_f \equiv z$ is not continuous.
Thus, it is not clear how to show that $V(U, p)$ is nonempty.
We are not aware of any simple proof of this fact in the multivariate case.
}
\end{remark}

\subsection{H\"older regularity in $\mathbf{L_p(\re^s)}$} \label{ss.p.Hoelder}

\noindent To be able to determine the exact H\"older regularity of a refinable function $\varphi \in L_p(\re^s)$,
$1 \le p <\infty$, we need to adjust the definitions of the transition matrices $T_d$, $d \in D(M)$.
To do so we replace the set $\Omega$ in Definition \ref{d.omega} by the set $\tilde{\Omega}$ in \eqref{eq.omega0}, the latter contains $\Omega$ and
is determined by a certain admissible absorbing set $\Delta$.

\begin{defi}\label{d.absorb}
Let $1 \le p <\infty$. A set $\Delta \subset \re^s$ is called absorbing if, for all $f \in L_p(\re^s)$,
the H\"older exponent of $f$ along $\Delta$ satisfies $\alpha_{f, \Delta , p} \, = \, \alpha_{f, p}$.
\end{defi}

\begin{remark}\label{r.280}{\em An arbitrary set that contains some neighborhood of the origin is absorbing.
It is also easy to show that any convex body (convex set with a nonempty interior) that contains
the origin is absorbing.}
\end{remark}

\noindent For the sake of simplicity, we choose $\Delta$ to be an arbitrary
simplex with one of the vertices at the origin and such that its interior intersects all the spaces $J_i, i = 1, \ldots , q(M)$. In this case $\Delta$ is absorbing, and the sets $\Delta \cap J_i$, $i=1,\ldots,q(M)$, are absorbing in the corresponding subspaces $J_i$. We call such a simplex $\Delta$ {\em admissible}. Note that for each $t > 0$,
the set $t\, \Delta$ is also an admissible simplex.

\noindent Define $\tilde{\Omega} \subset \z^s$ to be the minimal set such that
\begin{equation}\label{eq.omega0}
K_{\Gamma}\, + \, \Delta \  \subset \ \bigcup\limits_{k \in \tilde{\Omega}} (k + G)\, .
\end{equation}
Such a set $\, \tilde{\Omega}\, $ always exists, due to  $\displaystyle \cup_{k \in \z^s} (k + G) \, = \, \re^s$. Note that $\Omega \subseteq \tilde \Omega$. In many cases $\tilde \Omega = \Omega$, but not always,
see examples~\ref{ex.p10} and \ref{ex.p20}.
 Thus,
 $\,  {\rm supp} \, \varphi \, + \, \Delta \  \subset \  \tilde{\Omega}+G\, .$ Using $\tilde{\Omega}$ we redefine
$$
\tilde{v}_{\varphi} \, = \, \bigl(\varphi(\cdot +k)\bigr)_{k \in \tilde{\Omega}} \quad \hbox{a.e. on} \quad \re^s,
$$
and
$$
 (\tilde{T}_d)_{a,b}=c_{Ma-b+d}, \quad a,b \in \tilde{\Omega}, \quad d \in D(M),
$$
are now of size  $\tilde{N} = |\tilde{\Omega}|$. This leads to the appropriate modification
$$
 \tilde{\cA}=\{ \tilde{T}_d|_U \ : \ d \in D(M)\}
$$
of the finite set $\cA$ in \eqref{def.cA}. The modified subspaces $\tilde{V}$, $\tilde{U}$, $\tilde{U_i}$, $i=1, \ldots, q(M)$, and $\tilde{W}$  differ from the subspaces $V, U, U_i$ and $W$, respectively, only by the lengths of their corresponding elements. We are now ready to formulate the main result of this subsection.

\begin{theorem}\label{th.holder-p}
For a refinable function~$\varphi \in L_p(\re^s)$, $1 \le p <\infty$,
\begin{equation}\label{eq.holder-direct-p}
\alpha_{\varphi, J_i , p}\quad = \quad  \log_{\, 1/r_i} \, \rho_{p}(\tilde{\cA}|_{\tilde{U}_i})\, , \qquad  i = 1, \ldots , q(M),
\end{equation}
and, consequently,
\begin{equation}\label{eq.holder-p}
\alpha_{\varphi , p}\quad = \quad   \min\limits_{i = 1, \ldots , q(M)}\,
\log_{\, 1/r_i} \, \rho_{p}(\tilde{\cA}|_{\tilde{U}_i})
\end{equation}
\end{theorem}
\begin{proof} We first show that $\alpha_{\varphi, J_i, p}
\ge \log_{\, 1/r_i} \, \rho_{p}(\tilde{\cA}|_{\tilde{U}_i})$. Set $\rho_{i,p}=\rho_{p}(\tilde{\cA}|_{\tilde{U}_i})$.
Choose an arbitrary $\tilde{h} \in J_i \cap \Delta$ such that
$\|\tilde{h}\| < \delta$, $\delta \in (0,1)$. Then the function $\psi = \varphi(\cdot+\tilde{h}) - \varphi$
is supported on $K_{\Gamma}+\Delta$. Hence, the vector-valued function $\tilde{v}_{\psi}$ is well defined
on~$G$. Thus, for arbitrary $k \in \n$, we have
\begin{eqnarray*}
 && \bigl\| \varphi(\cdot + M^{-k} \tilde{h})\, - \, \varphi \bigr\|_{L_p(\re^s)}^p\ = \
\bigl\| \bT^k \bigl( \varphi(\cdot + \tilde{h})\, - \, \varphi\bigr) \bigr\|_{L_p(\re^s)}^p\ = \
\bigl\| \bT^k  \psi \bigr\|_{L_p(K_{\Gamma}+\Delta)}^p\ =  \\
&& \hspace{0.3cm} \bigl\| \tilde{\bA}^k  \tilde{v}_{\psi} \bigr\|_{L_p(G)}^p\  =
\sum_{d_1, \ldots, d_k \in D^k(M)}\,  \int_{G_{d_1 \ldots d_k}}\| \tilde{A}_{d_1}\cdots \tilde{A}_{d_k} \tilde{v}_{\psi}
\bigl(M_{d_k}\cdots M_{d_1}x \bigr)\bigr\|^p\, dx\ = \\
&& \hspace{0.3cm}  \sum_{d_1, \ldots,d_k \in D(M)} \, m^{-k}\, \int_{G}\| \tilde{A}_{d_1}\cdots \tilde{A}_{d_k}
\tilde{v}_{\psi}(y)\bigr\|^p\, dy\ = \
\int_{G} \cF_p \, \bigl(k\, , \, \tilde{v}_{\psi}(y ) \, \bigr)\, dy\, .
\end{eqnarray*}
By Corollary~\ref{c.p-rad}, for every $\varepsilon>0$, there exists a constant $C(\varepsilon)>0$ such that
$$
\int_{G} \cF_p \, \bigl(k\, , \, \tilde{v}_{\psi}(y ) \, \bigr)\, dy \ \le \
C(\varepsilon) \, (\rho_{ i, p} + \varepsilon)^{kp} \|\tilde{v}_{\psi}\|^p_{L_p(G)}\ \le \
C(\varepsilon) \, (\rho_{ i, p} + \varepsilon)^{kp} \,
2^p\, \|\tilde{v}_{\varphi}\|^p_{L_p(G)} \, .
$$
Thus, we obtain
\begin{equation} \label{aux70}
\bigl\| \varphi(\cdot + M^{-k} \tilde{h})\, - \, \varphi \bigr\|_{L_p(\re^s)} \ \le \ \tilde{C}(\varepsilon)\, (\rho_{i, p} + \varepsilon)^k, \quad k \in \n,
\end{equation}
where $\tilde{C}(\varepsilon)>0$ is independent of either $k$ or $\tilde{h}$. Choose an arbitrary $h \in J_i \cap \Delta , \, \|h\| < \delta$,
and let $k \in \n$ be the largest integer such that $\|M^kh\| < \delta$. From $\|M^{k+1}h\| \ge \delta$ and substituting $\tilde{h} = M^k h$ in \eqref{aux70}, we obtain
$$
\bigl\| \varphi(\cdot + h)\, - \, \varphi \bigr\|_{L_p(\re^s)} \ \le \ \tilde{C}(\varepsilon)\, (\rho_{i, p} + \varepsilon)^k\, ,
\qquad k \ge C_1+ \log_{\, r_i} (\delta /h)\,
$$
for some constant $C_1>0$.  Combining these estimates, we obtain $\alpha_{\varphi, J_i, p}\, \ge \, \log_{\, 1/r_i}\,
 (\rho_{i, p} + \varepsilon)$. Taking the limit for $\varepsilon \to 0$, we obtain the claim.

To establish the reverse inequality~$\alpha_{\varphi, J_i, p}
\le \log_{\, 1/r_i} \, \rho_{i, p}$, we argue as in the second part of the proof of Theorem~\ref{th.p}. We show the existence of a vector $h \in J_i \cap \Delta$ and of a subset $H \subset J_i$ of positive Lebesque measure (on the space~$J_i$) for which inequality~(\ref{eq.p2}) holds
(with $\rho_p$ replaced by $\rho_{i, p}$). Taking a limit in that inequality as $k \to \infty$
and using the fact that $k \, \le \, C_2 + \log_{\, r_i}\, \frac{\|h\|}{\|h_k\|}$, where
$C_2>0$ independent of $k$, we complete the proof.
\end{proof}

\subsection{Examples}\label{ss.p.ex}

\noindent The following examples illustrate  the need for the modifications of the set $\Omega$ in
subsection~\ref{ss.p.Hoelder}.

\begin{ex}\label{ex.p10}
{\em The solution of  the univariate refinement equation
$$
\varphi(x) = \varphi(2x) + \varphi(2x-1), \quad x \in \re,$$
 is the characteristic function $\varphi=\chi_{[0,1)}$. The $L_p$-regularity of $\varphi$ is
$\alpha_{\varphi, p} = \frac1p$. In this case, $M=2$  and, for the standard set of dyadic digits $D(M) = \{0, 1\}$, we have $G = [0,1]$ and $\Omega = \{0\}$. 
Hence, $N=1$ and we get $T_0 = T_1 = 1$. The common invariant subspace of $T_0$ and $T_1$ is trivial $U = \{0\}$, hence,
by definition, $\rho_p(\cA) = 0$. Thus, $\alpha_{\varphi, p} = \frac{1}{p}$, while
$\log_{1/2} \rho_p(\cA)\, = \, +\infty$. We see that $\alpha_{\varphi, p} \ne \log_{1/2} \rho_p(\cA)$.

\noindent On the other hand, for $\tilde{\Omega} = \{0, 1\}$, we get
$$
\tilde{T}_{0} \ = \ \left(
\begin{array}{cc}
1 & 0 \\
0 & 1
\end{array}
\right) \qquad \hbox{and} \qquad
\tilde{T}_{1} \ = \ \left(
\begin{array}{cc}
1 & 1 \\
0 & 0
\end{array}
\right).
$$
The corresponding common invariant subspace $\tilde{U} = \tilde{W} = \{u \in \re^2 \ : \ u_1 + u_2 = 0\}$ is one dimensional, and
$\tilde{A}_0 =  \tilde{A}_1 = 1$. Clearly, $\rho_p (\{\tilde{A}_0, \tilde{A}_1\}) = 2^{-1/p}$. By Theorem~\ref{th.holder-p},
we obtain the correct H\"older exponent $\alpha_{\varphi , p} = \frac1p$.
}
\end{ex}

\noindent The next example shows that in some cases $\Omega=\tilde{\Omega}$.

\begin{ex}\label{ex.p20}
{\em The solution of the univariate refinement equation
$$
\varphi(x) \ = \ \varphi(3x) \ + \ \varphi(3x-1)\ + \ \varphi(3x-5), \quad x \in \re,
$$
is $\varphi=\chi_{G}$, where $G$ is the tile in~$\re$ corresponding to the dilation $M = 3$ and to the digit set
$D(M) = \Gamma = \{0, 1, 5\}$. Thus, $\hbox{supp}(\varphi) \subseteq [0,\frac{5}{2}]$.

\noindent For the standard  set of triadic digits $D(M) = \{0, 1, 2\}$, we have
$G = [0,1]$ and  $\Omega = \{0, 1, 2\}$, i.e. $N=3$.
In this case, $K+\Delta \subseteq [0, 3]$. Hence, one can take
the complemented set $\tilde{\Omega}=\Omega = \{0,1,2\}$. 

}
\end{ex}

\begin{remark}\label{r.p100}
{\em It is well known that, if $p$ is an even integer, then $\rho_p(\cA)$
can be efficiently computed as an eigenvalue of a certain matrix derived from the
matrices in~$\cA$ \cite{JZ,P97}. Hence, Theorem~\ref{th.holder-p} allows us to find
the H\"older $L_p$-regularity at least for even integers $p$, in particular, for $p=2$, see
Example \ref{p.ex30}. }
\end{remark}

\begin{ex}\label{p.ex30}
{\em For the refinement equation from subsection~\ref{ss.examples.C}, we have $\tilde \Omega = \Omega$ and
$U = W$.  Therefore, Theorem~\ref{th.holder-p} yields $\alpha_{\varphi, p}\, = \,
\log_{1/\rho(M)} \rho_p(\cT|_{W})$. Furthermore, $\rho(M) = \frac{1+\sqrt{13}}{2} \, = \, 2.30277\ldots$
and $\rho_2 (\cT|_{W})\, = \, 0.79736\ldots$. Hence, $\alpha_{\varphi, 2}\, = \, \log_{1/\rho(M)} \rho_2(\cT|_{W})\, = \, 0.27148\ldots$. Recall  that the H\"older exponent in $ C(\re^s)$ of $\varphi$
is $\alpha_{\varphi}\, = \, 0.07652 \ldots$.  }
\end{ex}



\medskip

\noindent \textbf{Acknowledgements.} The authors are grateful to N.~Guglielmi and T.~Mejstrik who kindly spent their valuable time to
help us with computation issues.



\begin{thebibliography}{}

\bibitem{Hyp_W}
P.~Abry, M.~Clausel, S.~Jaffard, S.~Roux, B.~Vedel,
\newblock {\em Self-similar anisotropic texture analysis: the hyperbolic Wavelet transform contribution},
\newblock IEEE Trans. Image Process, 22 (2013), 4353 - 4363.

\bibitem{B}
N.\,E.~Barabanov,
\newblock {\em Lyapunov indicator for discrete inclusions, I-III},
\newblock Autom. Remote Control, 49 (1988), 152 - 157.

\bibitem{Bownik}
M.~Bownik,
\newblock{\em Anisotropic Hardy spaces and wavelets},
\newblock Memoirs Amer. Math. Soc. 164 (2003), no.~781.


\bibitem{CHM2}
C.\,A.~Cabrelli, C.~Heil, U.\,M.~Molter,
\newblock {\em Self-similarity and multiwavelets in higher dimensions},
\newblock Memoirs Amer. Math. Soc., 170 (2004), no.~807.
\smallskip

\bibitem{CDM}
A.~S.\,Cavaretta,  W.\,Dahmen,  C.A.\,Micchelli,
\newblock  {\em Stationary subdivision},
\newblock Memoirs Amer. Math. Soc. 93 (1991), no.~453.

\bibitem{Charina}
M. Charina,
\newblock {\em  Vector multivariate subdivision schemes:
Comparison of spectral methods for their regularity analysis},
\newblock Appl. Comp. Harm. Anal., 32 (2012), 86 - 108.

\bibitem{CharinaDRT}
M. Charina, M.~Donatelli, L.~Romani, V.~Turati,
\newblock {\em Multigrid methods: Grid transfer operators and subdivision schemes},
\newblock Linear Algebra Appl., 520 (2017), 151 - 190.

\bibitem{CharinaDRT2}
M. Charina, M.~Donatelli, L.~Romani, V.~Turati,
\newblock {\em Anisotropic, interpolatory subdivision and multigrid},
\newblock arXiv:1708.03469

\bibitem{CharinaProtasov}
M.~Charina, V.~Yu.~Protasov,
\newblock {\em Smoothness of anisotropic wavelets, frames and subdivision
schemes},
\newblock arXiv:1702.00269

\bibitem{CJR}
D.-R.\,Chen, R.-Q.\,Jia, S.D.Reimenschneider,
\newblock {\em Convergence of vector subdivision schemes in Sobolev spaces},
\newblock Appl. Comp. Harm. Anal., 12 (2002), 128 - 149.
\smallskip


\bibitem{Chr}
O.~Christensen,
\newblock {\em  An Introduction to Frames and Riesz Bases},
Birkh\"auser-Verlag, Basel, 2003.

\bibitem{Chui}
C.~Chui,
\newblock {\em  An Introduction to Wavelets},
Academic Press, London, 1992.

\bibitem{ChuiV}
C.~Chui, J.~de~Viles,
\newblock {\em  Wavelet Subdivision Methods: GEMS for Rendering Curves and Surfaces},
CRC Press, 2011.

\bibitem{CD}
A.~Cohen, I.~Daubechies,
 \newblock {\em A new technique to estimate the regularity of refinable functions},
  \newblock Revista Mathematica Iberoamericana, 12 (1996), 527 - 591.

\bibitem{CGV}
A. Cohen, K. Gr\"ochenig, L. Villemoes,
\newblock {\em Regularity of multivariate refinable functions},
\newblock Constr. Approx., 15 (1999), 241 - 255.

\bibitem{CoiW}
R.~R.~Coifman, G.~Weiss,
\newblock {\em Extensions of Hardy spaces and their use in analysis},
\newblock Bull. Am. Math. Soc., 83, 1997.

\bibitem{CH}
D.~Collela, C.~Heil,
\newblock {\em Characterization of scaling functions: I. Continuous solutions},
\newblock SIAM J. Matrix Anal. Appl., 15 (1994), 496 - 518.

\bibitem{CoSauer}
M.~Cotronei, D.~Ghisi, M.~Rossini, T.~Sauer,
\newblock {\em An anisotropic directional subdivision and multiresolution scheme},
\newblock Adv. Comput. Math., 41 (2015), 709 – 726.

\bibitem{Daub}
I.~Daubechies,
\newblock {\em Ten Lectures on Wavelets},
CBMS-NSF Regional Conference Series in Applied Mathematics, vol.~61,
SIAM, Philadelphia, 1992.


\bibitem{DL}
I.~Daubechies, J.~Lagarias,
\newblock {\em Two-scale difference equations.
II. Local regularity, infinite products of matrices and fractals},
\newblock SIAM J. Math. Anal., 23 (1992), 1031 - 1079.

\bibitem{DDL}
G.~Derfel, N.~Dyn, A.~Levin,
\newblock {\em Generalized refinement equations and subdivision processes},
\newblock J. Approx. Theory, 80 (1995),  272 - 297.

\bibitem{DD}
G.~Deslauriers, S.~Dubuc,
\newblock {\em Symmetric iterative interpolation processes},
\newblock Constr. Approx., 5 (1989), 49 - 68.

\bibitem{DH}
T.~B.~Dinsenbacher, D.~P.~Hardin,
\newblock{\em Nonhomogeneous refinement equations},
\newblock in Wavelets, Multiwavelets, and their Applications, A. Aldroubi and E. Lin, eds., AMS, Providence, RI,
1998,  117 - 127.

\bibitem{Dub}
S.~Dubuc,
\newblock {\em Interpolation through an iterative scheme},
\newblock J. Math. Anal. Appl.,  114 (1986),  185 - 204.

\bibitem{DynLevin02}
N.~Dyn, D.~Levin,
\newblock {\em Subdivision schemes in geometric modeling},
\newblock Acta Numer., 11 (2002), 73 - 144.


\bibitem{E}
T.\,Eirola,
\newblock {\em Sobolev characterization of solutions of dilation equations},
\newblock SIAM J. Math. Anal., 23 (1992), 1015 -1030.

\bibitem{FS}
D.-J.\,Feng, N.\,Sidorov,
\newblock {\em Growth rate for beta-expansions},
\newblock Monatsh. Math., 162 (2011), 41 - 60.



\bibitem{GroH}
K.~Gr\"ochenig, A.~Haas,
\newblock {\em Self-similar lattice tilings},
\newblock J. Fourier Anal. Appl., 2 (1994),  131 - 170.

\bibitem{GroM}
K.~Gr\"ochenig, W.R.~Madych,
\newblock {\em  Multiresolution analysis, Haar bases and self-similar tilings of $\re^n$},
\newblock  IEEE Trans. Inform. Theory, 38 (1992),  556 - 568.

\bibitem{GP1}
N.~Guglielmi, V.Yu.~Protasov,
\newblock {\em Exact computation of joint spectral characteristics of matrices},
\newblock Found. Comput. Math., 13 (2013),  37 - 97.

\bibitem{GP2}
N.~Guglielmi, V.Yu.~Protasov,
\newblock {\em Invariant polytopes of sets of matrices with applications to regularity of wavelets and subdivisions},
\newblock SIAM J. Matrix Anal. Appl., 37 (2016),  18 - 52.

\bibitem{HJ}
B.~Han, R-Q.~Jia,
\newblock {\em Multivariate refinement equations and convergence of subdivision schemes},
\newblock SIAM J. Math. Anal., 29 (1998), 1177 - 1199.

\bibitem{H1}
B.\,Han,
\newblock {\em Computing the smoothness exponent of a symmetric multivariate refinable function},
\newblock  SIAM J. Matr. Anal. Appl., 24 (2003), 693 - 714.

\bibitem{H2}
B.~Han,
\newblock {\em Vector cascade algorithms and refinable function vectors in Sobolev spaces},
\newblock J. Approx. Theory,  124 (2003), 44 - 88.

\bibitem{H3}
B.\,Han,
\newblock {\em Solutions in Sobolev spaces of vector refinement equations with a general dilation matrix},
\newblock Advances Comput. Math., 24 (2006), 375 - 403.

\bibitem{J95}
R.-Q.\,Jia,
\newblock {\em Subdivision schemes in $L_p$ spaces},
\newblock Advances Comp. Math., 3 (1995), 455 - 454.


\bibitem{J}
R.-Q.\,Jia,
\newblock {\em Characterization of smoothness of multivariate refinable functions in Sobolev spaces},
\newblock Trans. Amer. Math. Soc., 351 (1999), 4089 - 4112.

\bibitem{JJ}
R.-Q.\,Jia, Q. Jiang,
\newblock {\em Spectral analysis of transition operators and its applications to smoothness analysis of wavelets},
\newblock SIAM J. Matrix Anal. Appl., 24 (2003), 1071 - 1109.

\bibitem{JJS}
R.-Q.\,Jia, Q. Jiang, Z. Shen,
\newblock {\em Distributional solutions of nonhomogeneous discrete and continuous refinement equations},
\newblock SIAM J. Math. Anal., 32 (2000), 420 - 434.

\bibitem{JZ}
R.-Q.\,Jia, S.R.~Zhang,
\newblock {\em Spectral properties of the transition operator associated with multivariate refinement equation},
\newblock Lin. Alg. Appl., 292 (1999), 155 - 178.

\bibitem{JiangLi}
Q.~Jiang,  B.~Li,
\newblock {\em Quad/triangle subdivision, nonhomogeneous refinement equation and polynomial reproduction},
\newblock Math. Comput. Simul., 82 (2012), 2215 - 2237.

\bibitem{JiangLiZhu}
Q.~Jiang,  B.~Li, W.~Zhu,
\newblock {\em Interpolatory quad/triangle subdivision schemes for surface design},
\newblock Comput. Aided Geom. Des., 26 (2009), 904 - 922.

\bibitem{KM}
R.\,Kapica, J.\,Morawiec,
\newblock {\em Refinement type equations and Grincevicjus series},
\newblock J. Math. Anal. Appl., 350 (2009), 393 - 400.

\bibitem{Prot_book}
A.~Krivoshein, V.~Yu.~Protasov, M.~Skopina,
\newblock {\em  Multivariate wavelet frames},
Springer, 2016.

\bibitem{LW1}
J.\,Lagarias, Y.\,Wang,
\newblock {\em Integral self-affine tiles in $\re^n$. II. Lattice tilings, }
\newblock J. Fourier Anal. Appl., 3 (1997), 83 - 102.

\bibitem{LW2}
J.\,Lagarias, Y.\,Wang,
\newblock {\em Corrigendum and addendum to: Haar bases for $L_2(\re^n)$ and algebraic
number theory},
\newblock J. Number Theory,  76 (1999), 330 - 336.

\bibitem{MR}
C.~M\"oller, U.~Reif,
\newblock {\em A tree-based approach to joint spectral radius determination},
\newblock Lin. Alg. Appl., 563 (2014), 154 - 170.

  
\bibitem{NPS}
I.\,Y.~Novikov, V.\,Yu.~Protasov,  M.\,A.~Skopina,
\newblock {\em Wavelets theory},
\newblock  AMS, Translations Mathematical Monographs, 239 (2011).

\bibitem{ReifPeter08}
J.~Peter, U.~Reif,
\newblock {\em Subdivision Surfaces, Geometry and Computing},
Springer-Verlag, Berlin, 2008.

\bibitem{P96}
V.\,Yu.~Protasov,
\newblock {\em The joint spectral radius and invariant sets of linear operators.}
\newblock Fundam. Prikl. Mat., 2 (1996),  205 - 231.

\bibitem{P97}
V.\,Yu.~Protasov,
\newblock {\em The generalized spectral radius. A geometric approach},
\newblock Izvestiya Math., 61 (1997), 995 - 1030.

\bibitem{P00}
V.\,Yu.~Protasov,
\newblock {\em Refinement equations with nonnegative coefficients},
\newblock J. Fourier Anal. Appl., 6 (2000),  55 - 78.

\bibitem{P06}
V.\,Yu.~Protasov,
\newblock {\em Fractal curves and wavelets},
\newblock Izvestiya Math., 70 (2006), 123 - 162.

\bibitem{P07}
V.\,Yu.~Protasov,
\newblock {\em Spectral decomposition of 2-block
Toeplitz matrices and refinement equations},
\newblock St.Petersburg Math. J., 18 (2007), 607 - 646.

\bibitem{P08}
V.\,Yu.~Protasov,
\newblock {\em Extremal $L_p$-norms and self-similar functions},
\newblock Lin. Alg. Appl.,  428 (2008),  2339 - 2357.

\bibitem{P12}
V.\,Yu.~Protasov,
\newblock {\em When do several linear operators share an invariant cone?},
\newblock Lin. Alg. Appl., 433   (2010), 781 - 789.

\bibitem{P16}
V.\,Yu.~Protasov,
\newblock {\em The Euler binary partition function and subdivision schemes},
\newblock Math. Comput., (2016), Published electronically: https://doi.org/10.1090/mcom/3128



\bibitem{Rio}
O.~Rioul,
\newblock {\em Simple regularity criteria for subdivision schemes},
\newblock  SIAM J. Math. Anal., 23 (1992), 1544 - 1576.

\bibitem{RonS}
A.~Ron, Z.~Shen,
\newblock {\em The Sobolev regularity of refinable functions},
\newblock  J. Approx. Theory, 106 (2000), 185 - 225.

\bibitem{RS}
G.\,C.~Rota, G.~Strang,
\newblock {\em A note on the joint spectral radius},
\newblock Kon. Nederl. Acad. Wet. Proc., 63 (1960),  379 - 381.

\bibitem{SauerK}
G.~Kutyniok, T.~Sauer,
\newblock {\em Adaptive directional subdivision schemes and shearlet multiresolution analysis},
\newblock SIAM J. Math. Anal. 41, (2009)  1436 - 1471.

\bibitem{Schmeisser_Triebel}
H.-J.~Schmeisser, H.~Triebel,
\newblock {\em Topics in Fourier analysis and functions spaces},
Leipzig: Akad. Verlagsges., 1987.


\bibitem{SZ}
G.~Strang, D.~X.~Zhou,
\newblock {\em Inhomogeneous refinement equations},
\newblock J. Fourier Anal. Appl., 4 (1998), 733 - 747.

\bibitem{Vil}
L.\,Villemoes,
\newblock {\em Wavelet analysis of refinement equations}
\newblock SIAM J. Math. Anal.,  25 (1994), 1433 - 1460.

\bibitem{W}
J.\,Warren, H. Weimer,
\newblock {\em Subdivision methods for geometric design},
Morgan-Kaufmann, 2002.

\end{thebibliography}
\end{document}